\documentclass{article}
\usepackage[utf8]{inputenc}
\usepackage[a4paper,left=3cm,right=1.3cm,top=3cm,bottom=4cm]
{geometry}
\usepackage{amsmath,amstext, amsthm, empheq, extpfeil,  
amsbsy,amssymb,marvosym,fancyhdr,graphicx,amscd,amsfonts,latexsym,delarray,stackrel,lineno,color,cite,appendix,xcolor,url,hyperref,slashed}
\usepackage{wasysym}
\usepackage{subfig}
\usepackage{endnotes}
\usepackage{enumitem}
\usepackage{srcltx}
\usepackage{mathrsfs}
\usepackage{float} 
\DeclareUnicodeCharacter{00D7}{\times}

\usepackage[most]{tcolorbox}
\tcbuselibrary{skins}

\newtcolorbox{graybox}{
  colback=gray!10,    % background color
  colframe=gray!80,   % frame color
  arc=2mm,            % rounded corners
  boxrule=0.5pt,      % border thickness
  left=6pt, right=6pt, top=6pt, bottom=6pt
}
\usepackage[all]{xy}
\usepackage{t1enc}
\usepackage{mathrsfs}
\usepackage{pifont}
\usepackage{yfonts}

\usepackage{calligra}
\usepackage{palatino}

\usepackage{mathpple}
\usepackage[T1]{fontenc}

\definecolor{Green}{rgb}{0,1,0}
\definecolor{Blue}{RGB}{0,0,191}
\definecolor{mathmodecolor}{RGB}{0,102,0}
\definecolor{keywordcolor}{RGB}{0,51,151}
\definecolor{sourcebackgroundcolor}{RGB}{255,247,223}
\definecolor{unixagred}{RGB}{255,0,0}
\definecolor{lightgray}{RGB}{191,191,191}
\definecolor{green}{RGB}{1,191,191}

\newcommand*\patchAmsMathEnvironmentForLineno[1]{%
  \expandafter\let\csname old#1\expandafter\endcsname\csname #1\endcsname
  \expandafter\let\csname oldend#1\expandafter\endcsname\csname end#1\endcsname
  \renewenvironment{#1}%
     {\linenomath\csname old#1\endcsname}%
     {\csname oldend#1\endcsname\endlinenomath}}%
\newcommand*\patchBothAmsMathEnvironmentsForLineno[1]{%
  \patchAmsMathEnvironmentForLineno{#1}%
  \patchAmsMathEnvironmentForLineno{#1*}}%
\AtBeginDocument{%
\patchBothAmsMathEnvironmentsForLineno{equation}%
\patchBothAmsMathEnvironmentsForLineno{align}%
\patchBothAmsMathEnvironmentsForLineno{flalign}%
\patchBothAmsMathEnvironmentsForLineno{alignat}%
\patchBothAmsMathEnvironmentsForLineno{gather}%
\patchBothAmsMathEnvironmentsForLineno{multline}%
}
\newtheorem{thm}{Theorem}[section]

\newtheorem{fact}[thm]{Fact}

\def\cS{\mathcal{S}}
\def\sheaf{{\mathscr{O} }}
\newcommand{\ie}{{\it i.e.\/}\ }
\newcommand{\eg}{{\it e.g.\/}\ }

\newcommand{\dln}{{D_{\log}^{(\lambda,N)}}}

\def\C{\mathbb C}

\def\sin{{{\rm sin}}}
\def\cos{{{\rm cos}}}

\def\dim{{\mbox{dim}}}
 \def\scal2{{\mathscr S}}

\def\sr0{{\cS^{\rm ev}_0}}

\def\sar0{{\cS_0(\A_\Q)}}

\newcommand{\GL}{{\rm GL}}
\def\pw{\operatorname{PW}}

 \def\sr0{{\cS^{\rm ev}_0}}
\def\A{{\mathbb A}}
\def\F{{\mathbb F}}

\def\Q{{\mathbb Q}}
\def\R{{\mathbb R}}

\def\Z{{\mathbb Z}}

\def\Tr{{\rm Tr}}

\def\cD{{\mathcal D}}
\def\cE{{\mathcal E}}
\def\cS{{\mathcal S}}

\def\Spec{{\rm Spec\,}}

\def\qqq{\,,\quad \forall}

\newcommand{\nil}[1]{}

\title{The Riemann Hypothesis: Past, Present and a Letter Through Time}
\author{[Alain Connes] }

\date{\today}

\begin{document}

\maketitle

\begin{abstract}
This paper, commissioned as a survey of the Riemann Hypothesis, provides a comprehensive overview of 165 years of mathematical approaches to this fundamental problem, while introducing a new perspective that emerged during its preparation. 
The paper begins with a detailed description of what we know about the Riemann zeta function and its zeros, followed by an extensive survey of mathematical theories developed in pursuit of RH—from classical analytic approaches to modern geometric and physical methods. We also discuss several equivalent formulations of the hypothesis.

Within this survey framework, we present an original contribution in the form of a "Letter to Riemann," using only mathematics available in his time. This letter reveals a method inspired by Riemann's own approach to the conformal mapping theorem: by extremizing a quadratic form (restriction of Weil's quadratic form in modern language), we obtain remarkable approximations to the zeros of zeta. Using only primes less than 13, this optimization procedure yields approximations to the first 50 zeros with accuracies ranging from $2.6 × 10^{-55}$ to $10^{-3}$. Moreover we prove a general result that these approximating values lie exactly on the critical line $\Re(z)=\frac 12$.

Following the letter, we explain the underlying mathematics in modern terms, including the description of a deep connection  of the Weil quadratic form  with the world of information theory. The final sections develop a geometric perspective using trace formulas, outlining a potential proof strategy based on establishing convergence of zeros from finite to infinite Euler products. While completing the commissioned survey, these new results suggest a promising direction for future research on Riemann's conjecture.
\end{abstract}
\tableofcontents

\section{Introduction}

The Riemann Hypothesis (RH) stands as perhaps the most famous unsolved problem in mathematics, having resisted all attempts at proof since Bernhard Riemann casually remarked in his 1859 memoir that "sehr wahrscheinlich" (very probably) all nontrivial zeros of his zeta function lie on the critical line $\Re(s) = 1/2$. 

\medskip

The essence of this paper is a three pages letter to Riemann, written as if it were possible to communicate across the years, and in which I will present a strategy (towards a proof of RH) which only uses mathematics that he was familiar with, together with the great power of modern computers.

\medskip

   In the section \textit{"Encounter with the Riemann zeta function"} we aim to convey to the reader a vivid picture of the extraordinary connection that Riemann unveiled between the complex zeros of the zeta function and the distribution of prime numbers. We begin with a short historical account of the prime number theorem, followed by Riemann’s explicit formula for $\pi(x)$, the number of primes below $x$, and its later generalizations into the  explicit formulas. From there we turn to one of the most striking consequences: Littlewood’s discovery that the difference $\pi(x)-\mathrm{Li}(x)$ changes sign infinitely often, an observation that overturned earlier expectations of a simple monotone law.\newline
We then trace the story of what has been established about the zeros of $\zeta(s)$, beginning with Hardy’s proof of their infinitude on the critical line, through Selberg’s breakthrough in 1942, and up to modern results concerning the proportion of zeros lying exactly on the critical line. To illuminate this theory, we also draw on the rich framework of entire functions, including Hadamard’s factorization theorem, Nevanlinna theory, and the remarkable Borchsenius--Jessen theorem, itself a success of the Bohr--Landau program on almost periodic functions.\newline
This part culminates with a result that would surely have amazed Riemann : Voronin’s universality theorem. It asserts that the zeta function, in the strip $1/2 < \Re(s) < 1$, has a chameleon-like property: by vertical translation it can approximate any preassigned function $f(s)$ that is continuous, non-vanishing, and holomorphic on a simply connected compact set within the strip.\newline
Taken together, these milestones show both the depth and the subtlety of the analytic structure of $\zeta(s)$. They prepare us to appreciate why the Riemann Hypothesis remains so resistant: the function itself embodies an astonishing range of behaviors, far more elusive than any first impression might suggest.\newline

   Riemann could certainly not have foreseen the extraordinary mathematical landscape that emerged from attempts to prove his conjecture. Entire theories were born from this pursuit with connections spanning from algebraic geometry to quantum physics. The  section \textit{"A century and a Half of Theory Building towards RH"} is devoted to a quick panorama of these theories.  The first is the extension of the problem into algebraic and arithmetic geometry. In this setting, the analogue of the Riemann zeta function in finite characteristic was developed, and the corresponding analogue of the Riemann Hypothesis was proved by Andr\'e Weil. To provide context, I will briefly recall the notions of Grothendieck schemes and \'etale cohomology. The essential point here is that Weil’s celebrated ``Rosetta stone'' analogy becomes clear when viewed through the lens of schemes. In the original Rosetta stone the parallel inscriptions are horizontal: the hieroglyphs at the top, the Demotic text in the middle, and the Greek at the bottom. Weil proposed a similar three-texts comparison for the Riemann Hypothesis. At the bottom, corresponding to the Greek text, lies Riemann’s work on Riemann surfaces. There is no analogue of the Riemann Hypothesis in this setting, but it illustrates the extraordinary power of transcendental methods in geometry. The middle text  corresponds to algebraic geometry over finite fields: here the analogue of the Riemann Hypothesis was established by Weil, with the zeta function arising as the generating function counting rational points on curves over finite field extensions. Finally, the top one, the ``hieroglyphic'' text, corresponds to the arithmetic case of the spectrum of $\mathbb{Z}$, still deeply mysterious.  

What is remarkable is that, thanks to Grothendieck’s theory of schemes, one sees a unifying framework across all three inscriptions: in each case one is dealing with regular schemes of dimension one. For this reason, I will briefly evoke schemes, \'etale cohomology, and motives.  

Another major generalization arises from the theory of automorphic forms and representation theory, where the Selberg trace formula and zeta functions play a central role.

I will then move on to random matrix theory and quantum chaos. A remarkable discovery, initiated by Montgomery and Dyson and later spectacularly confirmed by Odlyzko’s large-scale numerical experiments, is that the pair correlation of the normalized spacings between consecutive zeros of $\zeta(s)$ on the critical line coincides with the correlation of eigenvalues from the Gaussian Unitary Ensemble (GUE) of random matrix theory. Nicholas Katz and Peter Sarnak  extended this correspondence to entire families of L-functions and established a systematic theoretical framework for understanding their statistical behavior as explained in \S \ref{sectKS}.

 In making the comparison of zeta zeros with eigenvalues from the Gaussian Unitary Ensemble, one has to rescale locally to compensate for the varying density of zeros. This need for such an adjustment reveals the absence of a genuine ``ultraviolet model'' for the zeta function---a gap that will be addressed in the final part of this paper, in \S \ref{cmpro}.\newline  
A further advance came with the work of Keating and Snaith in 2000, who used random matrix theory to propose a conjectural formula for the moments of the Riemann zeta function.  

 I will then briefly turn to my own work from 1998 on the trace formula. The new ingredient there is that, instead of studying the zeta function directly, one considers the ideal it generates, which is equivalent to focusing on its zeros. Crucially, one can capture this ideal and the associated zeros without defining the function explicitly, nor appealing to its analytic continuation. The key is the construction of a geometric space---the ad\`ele class space---which acquires a central significance as the  quotient of the adèles by the ergodic multiplicative action of the rational numbers. This quotient reveals the zeros of $L$-functions, not only of the Riemann zeta function, as an absorption spectrum. Recent results underline a precise correspondence between the ad\`ele class space and the class field theory counterpart of the schemes intimately related to $\mathrm{Spec}\,\mathbb{Z}$ that Grothendieck’s theory of the \'etale site and the \'etale fundamental group unveiled as a generalization of Galois theory.  

The corresponding trace formula is analogous to Selberg’s, but more delicate, since it involves both an infrared and an ultraviolet cutoff. It is precisely here that, for the first time, prolate spheroidal wave functions naturally enter the scene.  

I will mention various generalizations of the Riemann zeta function, which broaden the perspective and suggest new avenues of exploration. For completeness, I will also touch upon $p$-adic functions: while no analogue of the Riemann Hypothesis is known in these cases, they have nonetheless given rise to a wealth of beautiful developments.  

Finally, I will turn in the section \textit{"Equivalent Formulations"} to the realm of equivalent formulations of the Riemann Hypothesis. The website/book  "Equivalences of the Riemann Hypothesis" by Kevin Broughan (2017, Cambridge University Press),  systematically catalogs over 100 equivalent formulations. This comprehensive work organizes equivalences by mathematical area:

\begin{itemize}
\item Elementary arithmetic functions
\item Prime counting formulations  
\item Analytic equivalences
\item Functional analysis criteria
\item Probabilistic statements
\item Matrix and operator theory
\item Dynamical systems
\end{itemize}
 Some of these equivalents of RH are strikingly simple : for instance, one lets 
  $R_n$ be the Redheffer matrix,  an $n \times n$ (0,1)-matrix with $R_{ij} = 1$ if $j = 1$ or $i$ divides $j$. Then 
  $$ RH \iff \det(R_n) = O(n^{1/2+\epsilon}), \ \ \forall \epsilon > 0$$
  
 These are alluring but treacherous: the history of the subject shows how easily an attempt to solve the problem can be drawn into the gravitational pull of this ``black hole'' of equivalences.\newline
 I will discuss two of these equivalent formulations, first Weil's positivity criterion, then Robin's criterion as improved by Lagarias. This latter criterion shows that RH belongs to the type of statements that Hilbert, in his program on the foundations of Mathematics, was hoping to be "provable if true". But Godel's theorem showed that precisely for this type of statement truth does not imply provability. This will be discussed in details involving Chaitin's theory of algorithmic complexity.

 The next part of this paper, entitled \textit{``A Letter to Professor Bernhard Riemann''}, 
is of a completely different nature. It is written in the form of a letter addressed 
to Riemann himself, in an imaginary dialogue across time. The purpose of the letter 
is to point out a simple but surprising observation: a strategy towards RH which 
fits entirely within Riemann’s own mathematical outlook, as exemplified for instance 
by his proof of the conformal mapping theorem. Importantly, this strategy requires 
only the mathematical concepts and tools available in 1859.

Concretely, it allows one to recover the first nontrivial zeros of the zeta function 
by using only a few factors from the Euler product. For example, truncating the 
Euler product at the prime $13$ and computing with the method described, one obtains 
an approximation of the first 50 zeros of $\zeta(s)$  whose precision is best understood by giving the probability that one could get such an approximation by chance. The likelihood 
that such agreement could occur merely \emph{by chance} is about $10^{-1235}$. 
To give a sense of scale, this is approximately the probability of correctly 
guessing the outcomes of over 4000 consecutive coin tosses --- a feat so 
improbable that, for all practical purposes, it excludes both coincidence 
and computational error.\newline
One might therefore be tempted to conclude: ``This is simply a new algorithm for 
computing the zeros of the Riemann zeta function.'' Were that the case, the Riemann 
Hypothesis itself would follow, since a general theorem ensures that whenever the 
smallest eigenvalue of the corresponding operator is simple and even\footnote{meaning that the associated eigenfunction is even}, 
the resulting approximating numbers form the spectrum of a self-adjoint operator 
and hence are all real. 
   \newline
This letter to Riemann is followed, in the section \textit{The Next Small Steps}, 
by a more detailed discussion of its content, together with an outline of a 
possible strategy to justify rigorously that, in general, the approximating 
zeros produced by this method converge to the actual zeros of the Riemann 
zeta function. While the numerical evidence in support of this is overwhelming, 
evidence alone is not a proof. In this natural line of attack, a central role is 
played by the prolate spheroidal wave functions, which I had introduced into 
this context in 1998. 

Computer experiments reveal a striking fact: when expressed in terms of 
spheroidal wave functions, one obtains an excellent approximation to the 
function that minimizes the Weil quadratic form. This observation provides 
strong motivation for the approach. It is of course possible that a complete 
proof along these lines may encounter serious obstacles. Yet regardless of 
how far one can progress, this direction naturally opens the way to a deeper 
exploration, whose starting point is presented in the section 
\textit{Geometric Perspectives} — namely, the unexpected relation 
between two seemingly distant mathematical worlds.\newline
On the one hand, there is the world of the Weil quadratic form. The central fact here is Weil’s remarkable equivalence: the Riemann Hypothesis is equivalent to the positivity of  certain quadratic forms which involve only finitely many primes. This is striking, since one might expect that addressing the Riemann Hypothesis would require controlling the full infinity of primes. Here, however, the problem reduces locally, to finite collections at a time. Moreover, in this setting, as mentioned earlier, one can exploit the general construction of functions whose zeros lie entirely on the critical line.
 \newline 
 On the other hand, there is the world of the prolate wave functions developed by David Slepian and collaborators, with roots in Claude Shannon’s work in communication theory. The key link, established in the seminal work of Slepian, Landau, and Pollak, is a classical second-order differential operator on the line: the prolate operator, obtained as a confluence from the Heun equation—something entirely familiar to Riemann’s mathematical universe. This operator plays a dual role. In the infrared regime, it provides a means of approximating the minimal eigenvector of the Weil quadratic form. At the other extreme, as explained in the final section of this paper, it also supplies a model for a self-adjoint operator whose spectrum reflects the ultraviolet behavior of the zeros of the Riemann zeta function.
 \newline 
 Taken together, these observations underline the importance of understanding situations where only finitely many primes are involved, and of clarifying the interplay between the Weil quadratic form and the prolate operator. This connection is established through the trace formula, which I proved in 1998. In the last section, I will reformulate this trace formula using precisely the same ingredients that arise in communication theory and in the work of Slepian and collaborators, thereby making explicit the bridge it provides between the Weil quadratic form and the prolate functions.
  \newline    
  Behind the scene there lies a geometric space of fundamental importance: 
the \emph{semilocal adele class space} $Y_S$, associated with a finite set 
$S$ of places of $\Q$ that includes the archimedean place. By construction, 
$Y_S$ is obtained as the quotient  
\[
   \prod_{v \in S} \Q_v/\Gamma, \ \ \Gamma:=\{\pm \prod p_v^{n_v}\mid n_v\in \Z\}
\]
of the finite product of local fields by the diagonal action of the group $\Gamma\subset \Q^*$ generated by the primes in $S$. 
Weil's positivity criterion shows that proving the required positivity 
in this semilocal setting is equivalent to establishing the Riemann 
Hypothesis (in fact, one also captures the case of $L$-functions with 
Grössencharacters). 

A key advantage of this semilocal space, when compared with the full adele 
class space, is measure-theoretic: here the multiplicative and additive Haar 
measures are no longer singular with respect to one another, so the quotient 
behaves well from the viewpoint of measure theory. Topologically, however, 
the situation is more subtle. For each prime $p \in S$ one encounters a 
periodic orbit of length $\log p$, and these orbits encode precisely the 
contribution of $p$ to the explicit formula. Finally in \S \ref{cmpro}, I explain the recent discovery in our joint work with H. Moscovici \cite{CM}, where we show that the eigenvalues of the selfadjoint extension  of the prolate spheroidal operator reproduce the ultraviolet behavior of the squares of zeros of the Riemann zeta function and construct an isospectral family of Dirac operators whose spectra have the same ultraviolet behavior as the zeros of zeta.

\section{Encounter with the Riemann zeta function}

\subsection{Classical Analytic Number Theory}

    \subsubsection{Prime Number Theorem (PNT):}\label{sectpnt} In 1852, Pafnuty  Chebyshev, proved in \cite{Chebyshev} that the number $\pi(x)$ of primes less than $x$ fulfills 
    $$A \frac{x}{\log (x)}<\pi(x)<\frac{6 A}{5} \frac{x}{\log (x)}, \ \ \  A \approx 0.92129$$ for $x$ large enough, which allowed him to prove the convergence of the series  over the  primes $p$, $$\frac{1}{2 \log 2}+\frac{1}{3 \log 3}+\frac{1}{5 \log 5}+\frac{1}{7 \log 7}+\ldots$$
    He introduced two key counting functions $\vartheta$ and $\psi$, given as the sums 
    $$
    \vartheta(x):=\sum_{p<x} \log p, \  \  \psi(x):=\vartheta(x)+\vartheta(x^{1/2})+\vartheta(x^{1/3})+\ldots
    $$
     proved the  identity 
    $
    \sum_1^x \psi(x/n)=\log(x!)$
        and then used Stirling's formula. 
     The fact that the asymptotic relation $\pi(x) \sim x / \log x$ is equivalent to $\vartheta(x) \sim x$ was surely known to Chebyshev and was part of the math folklore  at that time. In fact in their  papers \cite{Hadamard}, \cite{ValeePoussin}, both Hadamard and de la Vallée Poussin prove $\vartheta(x) \sim x$ and do not mention explicitly the equivalent result $\pi(x) \sim x / \log x$.  Their independent proofs crucially established the non-vanishing of $\zeta(s)$ on the line $\mathrm{Re}(s) = 1$, using complex analysis techniques including Hadamard's theory of entire functions. The connection to RH is profound: while PNT requires only that $\zeta(s) \neq 0$ for $\Re(s) = 1$, the error term $\pi(x)-\operatorname{Li}(x)$  in the prime counting function is directly controlled by the location of the zeros. Under RH, the error term improves to $O(x^{1/2}\log x)$. The approximation by the integral logarithm differs substantially from its first term, one has
     $$\operatorname{Li}(x) \sim \frac{x}{\log x}+\frac{x}{\log ^2 x}+\frac{2 x}{\log ^3 x}+\cdots$$
     
     The proof of Hadamard that $\zeta(s) \neq 0$ for $\Re(s) = 1$ contains a beautiful idea that can be understood through the phases of prime powers and  reduced to $(-1)^2=1$. Indeed if for some real number $t$ one had $\zeta(1+it)=0$, the complex numbers $p^{-it}$ would tend to accumulate at $-1$, using the expansion of $-\log(1-p^{-it+\epsilon})$. But then  the complex numbers $p^{-2it}$ would tend to accumulate at $1$, thus creating a pole of the zeta function at $1+2it$ which is impossible unless $2t=0$.\newline
     The next step in the proof of the PNT by Hadamard and de la Vallée Poussin  stemmed from their need to extract precise asymptotic information about the distribution of primes from analytic properties of the zeta function. They  employed sophisticated techniques from the theory of entire functions and this required detailed estimates of the zeta function's behavior in various regions of the complex plane. While these methods were mathematically sound and represented brilliant mathematical insights, they obscured the conceptual relationship between the analytic input and the arithmetic output, \ie the asymptotic formula. 
     Edmund Landau provided the first major conceptual simplification by introducing Tauberian methods to prove the Prime Number Theorem, work that appeared in his influential 1909 "Handbuch der Lehre von der Verteilung der Primzahlen". Landau's innovation involved recognizing that the connection between analytic properties of the zeta function and arithmetic information about prime distribution could be mediated through a general class of theorems known as Tauberian results. These theorems provide conditions under which asymptotic properties of generating functions translate into asymptotic properties of their coefficient sequences. The progression from Landau's initial Tauberian approach through the Wiener-Ikehara theorem\endnote{The Wiener-Ikehara theorem states that if $A(x)$ is a non-negative, monotonically increasing function and $$f(s)=\int_0^\infty A(x) e^{-x s} d x$$ converges for $\operatorname{Re}(s)>1$, and if $f(s)-c /(s-1)$ has a continuous extension to $\operatorname{Re}(s) \geq 1$ for some constant $c \geq 0$, then $$\lim_{x \rightarrow \infty} e^{-x} A(x)=c.$$ This result provided exactly the framework needed to convert information about the zeta function's non-vanishing on $\operatorname{Re}(s)=1$ into precise asymptotic statements about prime-counting functions. 
The application to the Prime Number Theorem involves expressing the logarithmic derivative of the zeta function as a Mellin transform of the Chebyshev function $\psi(x)=\sum_{p^ k \leq x} \log p$ \ie the equality $-\zeta^{\prime}(s) / \zeta(s)=s \int_1^\infty \psi(x) x^{-(s+1)} d x$. The non-vanishing of $\zeta(s)$ on $\operatorname{Re}(s)=1$, combined with the known behavior of $\zeta(s)$ at $s=1$, ensures that $\zeta^{\prime}(s) / \zeta(s)+1 /(s-1)$ has the required analytic properties for applying the Wiener-Ikehara theorem. The theorem then directly implies that $\psi(x) \sim x$, from which the Prime Number Theorem follows through elementary arguments} to Newman's simplified presentation shows how mathematical understanding can deepen through the development of more general and conceptually transparent frameworks.\newline
      The 1949 elementary proof of the PNT by Selberg and Erdos showed that complex analysis is not logically necessary for proving the theorem.  The Selberg sieve, while developed for his elementary proof of PNT \cite{SelbergElementary}, became a fundamental tool in analytic number theory. The Selberg symmetry formula
    \[
   \sum_{p \leq x} \log ^2 p+\sum_{p q \leq x} \log p \log q=2 x \log x+O(x)
    \]
    exemplifies his ability to find unexpected simplicities in prime distribution.

    \subsubsection{Riemann's formula, von Mangoldt paper}\label{sectmangoldt}
    In his fundamental paper\footnote{which he sent to Chebyshev} \cite{Riemann}, Riemann was very careful to define which branch of the integral logarithm he was using and to state the conditional convergence\footnote{He was a master of conditional convergence and gave the first rigorous proof that conditionally convergent series can be rearranged to converge to any prescribed value or to diverge entirely. Riemann's rearrangement theorem emerged from his 1854 habilitation thesis "Über die Darstellbarkeit einer Function durch eine trigonometrische Reihe",  a work primarily focused on extending the theory of Fourier series to more general classes of functions. Published posthumously in 1867, \cite{Riemann1}, through Dedekind's editorial efforts.} in his formula for the function $f(x)$ obtained as the sum of the counting functions $\pi(x^{1/n})$
    \begin{equation}\label{piriemann0}
f(x)=\operatorname{Li}(x)-\sum_\alpha\left(\operatorname{Li}\left(x^{\frac{1}{2}+\alpha i}\right)+\operatorname{Li}\left(x^{\frac{1}{2}-\alpha i}\right)\right)+\int_x^{\infty} \frac{1}{t^2-1} \frac{d t}{t \log t}-\log 2
 \end{equation}
    where the $\frac{1}{2}+\alpha\ i$ are the non-trivial zeros of $\zeta$ with positive imaginary part and the order of the terms corresponds to increasing values of $\Re(\alpha)$.\newline
    As I will briefly mention at the beginning of my letter to Riemann, it is unfortunate that the care about the various branches of the integral logarithm has disappeared in some modern treatise such as the classic book of Edwards \cite{Edwards}, or in some technical papers, where \eqref{piriemann0} is written as above while there is no function $Li(z)$ or $Li(z)+\overline{Li(z)}$ for which the general term of the series tends to $0$. The problem is with the notation since if one follows Riemann closely one finds that the value he considers depends not just on $x^\rho$ but upon $\rho \log(x)$.      In fact in the paper of von Mangoldt \cite{vonMangoldt}, he correctly defines the meaning that he gives, not to the function $Li(z)$ but to the function $L i\left(e^w\right)$. Here is what he writes:\newline 
    \begin{small}
\begin{quote}Setting
$
w = u + iv
$
where  u  and  v  denote real numbers, and again understanding $h$  as a positive real variable, one only needs to make the following definitions:\newline
	1.	If  $v > 0$ , then

$$
L i\left(e^w\right) = \lim_{h \to \infty} \int_{-h+w}^w \frac{e^z}{z} dz + i\pi
$$
	2.	If  $v = 0$, then

$$
L i\left(e^w\right) = \lim_{h \to 0}\int_{-\infty}^{-h} \frac{e^z}{z} dz + \int_h^w \frac{e^z}{z} dz
$$
	3.	If  $v < 0$ , then

$$
L i\left(e^w\right) = \lim_{h \to \infty} \int_{-h+w}^w \frac{e^z}{z} dz - i\pi
$$
It is easily verified that the function  $L i\left(e^w\right)$ defined in this way indeed possesses the desired properties.\end{quote} 
\end{small}
Thus von Mangoldt was fully aware of the defect of the notation $Li(x^\rho)$ which unfortunately remained used carelessly, even though most authors are (probably) aware of the problem. \newline
One can test Riemann's formula, reformulated   using the function $\operatorname{Ei}(\rho\log x)$ in place of $\operatorname{Li}(x^\rho)$ (see \eqref{piriemann1}) with the help of the computer and one finds a pretty good agreement by taking the sum over a few thousands zeros of zeta.

 \subsubsection{Explicit formulas} \label{sectexplicit}
 These relate prime-counting functions to the nontrivial zeros of $\zeta(s)$, making precise the influence of the zeros on arithmetic functions. Von Mangoldt \cite{vonMangoldt} rigorously proved  Riemann's original explicit formula  as deduced from the simpler equality 
    \[
    \psi(x) = x - \sum_{\rho} \frac{x^{\rho}}{\rho} - \log(2\pi) - \frac{1}{2}\log(1-x^{-2})
    \]
    where $\psi(x) = \sum_{p^k \leq x} \log p$ is the Chebyshev function, and the sum runs over nontrivial zeros $\rho$. As in Riemann's case the $\sum_{\rho}$ is taken as the limit of the partial sums for $-T\leq \Im(\rho)\leq T$ and fails to be absolutely convergent. The error term is of the order of $(\log T)^2/T$ (see \cite{MV2007} Theorem 12.5).\newline   Guinand \cite{Guinand} and Weil \cite{WeilExplicit}  later developed more general explicit formulas connecting test functions and their Fourier transforms to distributions of zeros, providing a harmonic analysis perspective on the relationship between primes and zeros. The usual formulation for the Riemann zeta function involves test functions $f:\R_+^*\to \C$ fulfilling suitable regularity conditions\footnote{they are continuous and with continuous derivative except  at finitely many points where both $f(x)$ and $f'(x)$ have at most a discontinuity of the first kind, and at which the value of $f(x)$ and $f'(x)$ is defined as the average of the right and left limits. Moreover one assumes  that for some $\delta>0$ one has 
$$
f(x)=O(x^\delta), \ \text{for} \ x\to 0+, \ \ f(x)=O(x^{-1-\delta}), \ \text{for} \ x\to +\infty.
$$  }, and their Mellin transform  defined as
\begin{equation}\label{mellin}
 \tilde f(s):=\int_0^\infty f(x)x^{s-1}dx
 \end{equation}
They take the form 
 \begin{equation}\label{bombieriexplicit}
 \sum_{\rho}\tilde f(\rho)=\tilde f(0)+\tilde f(1)-\sum_v {\mathcal W}_v(f),
 \end{equation}
 where the  sum on the left hand side is over all complex zeros $\rho$ of the Riemann zeta function, and is in general only conditionally convergent while the sum on the right hand side   runs over all rational places  $v$ of $\Q$. The non-archimedean distributions $\mathcal W_p$ are defined, using $f^\sharp(x):=x^{-1}f(x^{-1})$, as 
 \begin{equation}\label{bombieriexplicit1}
 {\mathcal W}_p(f):=(\log p)\sum_{m=1}^\infty\left(f(p^m)+f^\sharp(p^m)\right)
 \end{equation}
while the archimedean distribution is given by 
 \begin{equation}\label{bombieriexplicit2}
 {\mathcal W}_\R(f):=(\log 4\pi +\gamma)f(1)+\int_{1}^\infty\left(f(x)+f^\sharp(x)-\frac 2x f(1)\right)\frac{dx}{x-x^{-1}}.
 \end{equation}

    \subsubsection{Changes of sign of $\pi(x) - \mathrm{Li}(x)$}Littlewood's famous 1914 result \cite{LittlewoodZeroFree} that $\pi(x) - \mathrm{Li}(x)$ changes sign infinitely often (assuming RH) demonstrated that the error term in the PNT exhibits complex oscillatory behavior, challenging the empirical observation that $\pi(x) < \mathrm{Li}(x)$ for all computed values. The first effective bound for such a sign change was given by Skewes \cite{Skewes1933}, who showed (assuming RH) that a crossing occurs before
    $$
    B=10^{10^{10^{34}}}
    $$
      This enormous bound has been dramatically improved  by Lehman \cite{Lehman}, who introduced the key idea of deriving an explicit formula for $u e^{-u / 2}\left(\pi\left(e^u\right)-\operatorname{Li}\left(e^u\right)\right)$ averaged by a Gaussian kernel. It was further improved  by  Riele \cite{teRiele1987} (to approximately $6.69 \times 10^{370}$), Bays and Hudson \cite{BaysHudson2000} (to approximately $1.4 \times 10^{316}$), and most recently by Chao and Plymen \cite{ChaoPlymen2010} and Saouter,  Demichel and Trudgian\cite{SaouterDemichel2010}, \cite{SaouterTrudgianDemichel}.
    
     \subsubsection{Hardy and Littlewood}\label{Hardy-Littlewood} Hardy's 1914 proof \cite{Hardy} that infinitely many zeros lie on the critical line used the functional equation and careful analysis of $Z(t) = e^{i\theta(t)}\zeta(1/2 + it)$, where $\theta(t)$ is chosen to make $Z(t)$ real for $t$ real,
     $$
\theta(t)=\operatorname{Im}\left(\log \Gamma\left(\frac{1}{4}+i \frac{t}{2}\right)\right)-\frac{t}{2} \log \pi
$$ By studying sign changes of $Z(t)$, he showed infinitely many zeros exist with $\Re(s) = 1/2$.

The Hardy-Littlewood work \cite{HardyLittlewood} on conditional results showed that RH implies numerous consequences about prime gaps  and other additive problems.  Their approach built upon Hardy's earlier work but incorporated more refined analytical techniques for estimating oscillatory integrals and asymptotic expansions. Their subsequent work involved careful analysis of the approximate functional equation for the zeta function\footnote{Siegel's 1932 paper on Riemann's notebooks shows that while Hardy and Littlewood had independently rediscovered the principal term of this development in 1920 through their "approximate functional equation" using similar methods (the saddle-point method), Riemann's work contained additional insights}, which provides an explicit representation for $\zeta(\mathrm{s})$ that is valid throughout the critical strip. By studying the cancellations and oscillations in this representation, they were able to demonstrate that the number of sign changes in Hardy's $Z$-function, and hence the number of zeros of $\zeta(\mathrm{s})$ on the  critical line of imaginary part in $[0,T]$, grows at least linearly with $T$. It was still unclear if a positive fraction of the number $N(T)$ of non-trivial zeros with imaginary part in $[0,T]$, estimated by Riemann as $N(T)\sim \frac{1}{2\pi}T \log T$,  lie on the critical line.

    \subsubsection{Selberg's early contributions (1940s–1950s):}\label{Selberg} The breakthrough that resolved this fundamental question came in 1942 through Atle Selberg's pioneering work, which introduced the innovative technique of mollification and proved for the first time that a positive proportion of nontrivial zeros lie on the critical line. The essential idea in Selberg's proof\endnote{Selberg's mollifier was of the form
$$
\sum_{n \leq \xi} \frac{\alpha_n\left(1-\frac{\log n}{\log \xi}\right)}{n^s}
$$
where the $\alpha_n$ are the Dirichlet series coefficients of $1 / \sqrt{\zeta(s)}$ with $\alpha_1=1$. He needed an average of $|\zeta(1 / 2+i t)|^2$ times the fourth power of this mollifier.} is to multiply the series for $\zeta(s)$ by the square of the partial sums of the series for $\zeta(s)^{-1/2}$. Selberg's achievement was to demonstrate that the number $\mathrm{N}_0(\mathrm{~T})$ of zeros on the critical line satisfies $\mathrm{N}_0(\mathrm{~T}) \gg \mathrm{T} \log \mathrm{T}$, which immediately implies that the ratio $\mathrm{N}_0(\mathrm{~T}) / \mathrm{N}(\mathrm{T})$ is bounded away from zero, establishing that a positive fraction of all zeros reside on the critical line. Selberg's work \cite{SelbergZeta, SelbergElementary} revolutionized the study of the zeta function through several key innovations.  His work on the distribution of $\log|\zeta(1/2 + it)|$ showed it follows a Gaussian distribution with specific mean and variance, presaging connections to random matrix theory.

    \subsubsection{Proportion of zeros on the critical line:} As explained above, Hardy's 1914 theorem \cite{Hardy}  was quantified by Selberg (1942) \cite{SelbergZeta} who showed that a positive proportion of zeros have real part $1/2$. Levinson's breakthrough in 1974 \cite{Levinson} proved that at least $1/3$ of the nontrivial zeros lie on the critical line, using mollifier techniques that smooth the zeta function near the critical line. Conrey's 1989 improvement \cite{Conrey1989} to $2/5$ used more sophisticated mollifiers and asymptotic analysis. The 2011 result by Bui, Conrey, and Young \cite{BuiConreyYoung} achieving $41\%$ employed multiple Dirichlet polynomials and optimization techniques from the theory of extremal functions, representing the current state of the art in this approach. Pratt et al \cite{Pratt} have the current record of $41.7 \%$ based on Feng's mollifier \cite{Feng}.

 \subsubsection{Zero-free regions and zero-density estimates:} Classical zero-free regions, such as the one proved by de la Vallée Poussin \cite{ValeePoussin} showing $\zeta(s) \neq 0$ for $\Re(s) \geq 1 - c/\log(|\Im(s)|+2)$, provide explicit bounds on the distribution of primes in arithmetic progressions. Zero-density estimates, pioneered by Bohr and Landau, bound the number $N(\sigma,T)$ of zeros with $\Re(s) \geq \sigma$ and $|\Im(s)| \leq T$. The best known results show that $N(\sigma,T) \ll T^{c(1-\sigma)}$ for various values of $c$ depending on $\sigma$, with improvement $$N(\sigma, T) \ll T^{\frac{3(1-\sigma)}{2-\sigma}}(\log T)^5$$ by Ingham \cite{Ingham1}. Montgomery and others provided increasingly sharp bounds as $\sigma \to 1$. The Guth-Maynard zero density result \cite{Guth} is the first improvement in Ingham's theorem near $\sigma=3 / 4$ in 85 years.

        \subsubsection{The Lindelöf Hypothesis:} The Lindelöf Hypothesis states that $\vert \zeta(1/2 + it) \vert \ll t^{\epsilon}$ for any $\epsilon > 0$. If $\sigma$ is real, then $\mu(\sigma)$ is defined to be the infimum of all real numbers $\alpha$ such that $\vert\zeta(\sigma+i t)\vert=\mathrm{O}\left(t^\alpha\right)$. One has $\mu(\sigma)=0$ for $\sigma>1$, and the functional equation  implies that $\mu(\sigma)=\mu(1-\sigma)-\sigma+1 / 2$. The Phragmén-Lindelöf theorem implies that $\mu$ is a convex function. The Lindelöf hypothesis states $\mu(1 / 2)=0$, which together with the above properties of $\mu$ implies that $\mu(\sigma)$ is 0 for $\sigma \geq 1 / 2$ and $1 / 2-\sigma$ for $\sigma \leq 1 / 2$.  While weaker than RH, it would still have significant applications to moments of $L$-functions and subconvexity problems.
        
        % \subsubsection{Landau–Siegel zeros:} A Landau-Siegel zero would be a real zero  exceptionally close to $s=1$ of a Dirichlet $L$-function $L(s,\chi)$ for a quadratic character $\chi$. While no such zeros are known to exist, their hypothetical existence would imply irregularities in the distribution of primes in arithmetic progressions \cite{LandauSiegel, HeathBrownSiegel}. The effective bounds on such zeros, if they exist, are notoriously weak. 
%\end{itemize}

\subsection{Theory of Entire and Meromorphic Functions}

   \subsubsection{Hadamard factorization theorem}\label{hadfact} An entire function $f$ has order $\sigma$ if
$$\sigma = \limsup_{r \to \infty} \frac{\log \log M(r)}{\log r}$$
where $M(r) = \max_{|z|=r} |f(z)|$. The Hadamard product formula states that an entire function $f$ of order 1 can be written as
$$f(z) = e^{A + Bz} \prod_{\rho} \left(1 - \frac{z}{\rho}\right) e^{z/\rho}$$
where the product is over all zeros $\rho$ of $f$, and $A$, $B$ are constants. 
The completed zeta function 
\begin{equation}\label{smallxi}
	\xi(s) = \frac{1}{2}s(s-1)\pi^{-s/2}\Gamma(s/2)\zeta(s)
\end{equation}
 is an entire function of order one with the Hadamard product representation\footnote{The infinite product is defined here as the limit as $T\to \infty$ of the product over $\rho$ with $\vert \rho\vert<T$ } \cite{HadamardFactor, Boas, LevinBook}
 \[
    \xi(z) = \xi(0) \prod_{\rho} \left(1 - \frac{z}{\rho}\right)
    \]
Riemann considers the function (he was denoting it as $\xi(s)$ but, since Hardy, one uses $\Xi(s)$)
\begin{equation}\label{bigxi}
	\Xi(s) = \xi(1/2 + is)
\end{equation}
 The functional equation $\xi(s) = \xi(1-s)$ implies $\Xi(s) = \Xi(-s)$, thus the function $\Xi$ is even and its zeros come in pairs $\pm\alpha$. Moreover, $\Xi$ is entire of order 1,  for each pair of zeros $\pm\alpha$, the exponential factors $e^{s/\alpha}$ and $e^{-s/\alpha}$ cancel and the Hadamard factorization of $\Xi$ takes the particularly simple form
$$\Xi(s) = \Xi(0) \prod_{\Re(\alpha) > 0} \left(1 - \frac{s^2}{\alpha^2}\right)$$
where the product is over  zeros with positive real part only. The two product representations are the same\footnote{using the identity $1-\frac{\frac{1}{2}+i s}{\frac{1}{2}+i a}=\frac{1-\frac{s}{a}}{1-\frac{i}{2 a}}$} except for the constant factor and while $\xi(0)=\frac 12$, one has $$\Xi(0)=-\frac{\zeta \left(\frac{1}{2}\right) \Gamma \left(\frac{5}{4}\right)}{2 \sqrt[4]{\pi }}\sim 0.497121$$ This nuance is at the origin of the only "typo" in Riemann's formula for the number of primes less than $x$, where he used $\log \Xi(0)$ in place of $\log \xi(0)$.

   \subsubsection{Nevanlinna theory} This value-distribution theory \cite{Nevanlinna1, NevanlinnaCollected, Laine} studies the growth and value distribution of meromorphic functions through characteristic functions and deficiency relations. It is worth noting that the Nevanlinna characteristic is closely related to Jensen’s formula. 
For an analytic function $f$ with zeros $a_k$ in $|z|\leq r$, Jensen’s formula asserts that, provided $f(0)\neq 0$,
\[
\log|f(0)| + \sum_{k} \log \frac{r}{|a_k|}
   = \frac{1}{2\pi} \int_0^{2\pi} \log |f(re^{i\theta})|\,d\theta.
\]
This relation between the average boundary growth of $f$ and the distribution of its zeros 
is precisely what Nevanlinna theory systematizes for meromorphic functions. 
The characteristic $T(r,f) = m(r,f) + N(r,f),  r > 0,$ can thus be viewed as a refined version of the quantities 
appearing in Jensen’s formula, with the counting function $N(r,f)$ recording the poles 
and the proximity function $m(r,f)$ recording the boundary growth. For the zeta function, one can define the Nevanlinna characteristic $T(r,\zeta)$ and study its growth\endnote{For a meromorphic function $f$ in the complex plane, the Nevanlinna characteristic is defined by
\[
T(r,f) = m(r,f) + N(r,f), \qquad r > 0,
\]
where
\[
m(r,f) = \frac{1}{2\pi}\int_{0}^{2\pi} \log^+\!\big|f\!\left(re^{i\theta}\right)\big|\,d\theta
\]
is the proximity function, and
\[
N(r,f) = \int_{0}^{r} \frac{n(t,f)-n(0,f)}{t}\,dt + n(0,f)\log r
\]
is the integrated counting function of poles, with $n(t,f)$ denoting the number of poles of $f$ (with multiplicity) in $|z|\leq t$. For an entire function, $N(r,f)=0$ and hence $T(r,f)=m(r,f)$.

In the case of the Riemann zeta function $\zeta(s)$, which is meromorphic in $\mathbb{C}$ with a single pole at $s=1$, one can define $T(r,\zeta)$ in the same way. A classical result of Borchsenius and Jessen (1948) shows that
\[
T(r,\zeta) \sim \frac{\log r}{2\pi}\int_{-\pi}^{\pi} \log^+\!\big|\zeta\!\left(re^{i\theta}\right)\big|\,d\theta
\qquad (r\to\infty),
\]
and, more precisely, they obtained the asymptotic
\[
T(r,\zeta) \sim \frac{1}{\pi}\,r\log r \quad \text{as } r\to\infty.
\]

Thus the Nevanlinna characteristic of $\zeta(s)$ grows essentially like $r\log r$, which is of the same order as the growth of $\log M(r,\zeta)$, where $M(r,\zeta)=\max_{|s|=r}|\zeta(s)|$.}. The theory provides tools like the First and Second Main Theorems that relate the distribution of $a$-points (solutions to $f(z) = a$) to the growth of the function. While direct applications to prove RH have not succeeded, Nevanlinna theory offers insights into the exceptional nature of the value $0$ for $\zeta(s)$ as shown by the Borchsenius-Jessen theorem discussed below, and provides analogies with Picard-type theorems \cite{GoldbergOstrovskii}. Recent work has explored connections between Nevanlinna theory and the universality properties of $\zeta(s)$.
   
   \subsubsection{Average behavior of $\log |\zeta(s)|$ and zeros of $\zeta(s)-x$}  The Borchsenius-Jessen theorem \cite{BorchJessen} provides a striking contrast to RH. While all nontrivial zeros of $\zeta(s)$ are conjectured to lie on $\Re(s) = 1/2$, the zeros of $\zeta(s) - x$ for any $x \neq 0$ cluster densely near the critical line from the right. Specifically:
    \begin{itemize}
        \item For $x \neq 0$, the equation $\zeta(s) = x$ has infinitely many solutions in any strip $1/2 < \Re(s) < \sigma.$
        \item The density of solutions increases without bound as $\Re(s) \to 1/2^+$
        \item If $N_x(T)$ counts zeros of $\zeta(s)-x$ with $\Re(s) > 1/2$ and $0 < \Im(s) < T$, then $N_x(T)/T \to \infty$ as $T \to \infty$
    \end{itemize}
    This phenomenon highlights the exceptional nature of the value $x = 0$ and shows that RH is unstable under perturbations: adding any constant to $\zeta(s)$ destroys the alignment of zeros on the critical line. The result connects to the mean value theorems for $\log \zeta(s)$ in the critical strip and  to the ergodic properties of the flow $s \mapsto s + it$ on the value distribution of $\zeta(s)$. The Borchsenius-Jessen theorem can be viewed as the crowning achievement of the Bohr-Landau program on almost periodic functions and Dirichlet series. Where Bohr and Landau established the general framework for understanding value distribution of Dirichlet series through almost periodicity, Borchsenius and Jessen achieved the definitive result for the Riemann zeta function itself, proving  
that $\log \zeta(\sigma+i t)$ has a continuous limiting distribution\endnote{The theorem states that for fixed $\sigma>1 / 2$, the values 
$$
\frac{1}{T} \int_0^T f(\log \zeta(\sigma+i t)) d t
$$
converge as $T \rightarrow \infty$ for continuous bounded functions $f: \mathbb{C} \rightarrow \mathbb{C}$
The limiting distribution $\mu_\sigma$ is a probability measure on $\mathbb{C}$ such that:\newline
- It has no atoms (continuous distribution)\newline 
- Its support is the whole complex plane when $1 / 2<\sigma<1$.\newline
Moreover in spite of the ambiguity in defining $\log \zeta(\sigma+i t)$ in the presence of zeros of zeta, the Borchsenius-Jessen theorem does not assume the Riemann Hypothesis.}.

\subsubsection{Voronin's Universality Theorem, zeta as a chameleon} This remarkable 1975 theorem \cite{Voronin} states that the Riemann zeta function can approximate any non-vanishing holomorphic function arbitrarily well in the critical strip. Voronin's original theorem  established universality for disks: If $0<r<1 / 4, f(s)$ is continuous on $|s| \leq r$, holomorphic on $|s|<r$ with $f(s) \neq 0$ on $|s| \leq r$, and $\epsilon>0$, then there exists $t_0$ such that
$$
\max _{|s| \leq r}\left|f(s)-\zeta\left(s+\frac{3}{4}+i t_0\right)\right|<\epsilon
$$
 Following Voronin's pioneering work, it was shown by Reich and Bagchi that if $K$ is a compact subset of the strip $\{s : 1/2 < \Re(s) < 1\}$ with connected complement, $f$ is continuous on $K$ and holomorphic in the interior of $K$ with $f(s) \neq 0$ on $K$, and $\epsilon > 0$, then there exists $t_0$ such that
    \[
    \max_{s \in K} |\zeta(s + it_0) - f(s)| < \epsilon
    \]
    Moreover, the set of such $t_0$ has positive lower density. Thus, the Riemann zeta function in the critical strip $\frac 12<\Re(z)<1$
 behaves like a mathematical chameleon: it can approximate any non-vanishing holomorphic function on compact sets simply by shifting vertically in the complex plane. This result uses, as in the Bohr-Landau almost periodic set-up, the ergodic properties  of the Kronecker flow  on the infinite dimensional torus which is the compact dual of the discrete group $\Q^*$ of non-zero rational numbers, but we now discuss the key ingredient which is of a different nature. \newline As we briefly mentioned above in a footnote, Riemann in his memoir on trigonometric series showed, in his discussion of the Dirichlet kernel, that a conditionally convergent series of real numbers can be given any sum $C$ after a suitable rearangement of its terms. What is amazing is that it is a generalization of this fact which is at the heart of Voronin's proof. 
The key result is the following (\cite{Laurincikas}, Theorem 1.16)
\begin{graybox}
Let $H$ be a Hilbert space, $\left\{x_m\right\}$ be a sequence in $H$ satisfying the two conditions:
 $$\sum_{m=1}^{\infty}\left\|x_m\right\|^2<\infty;
 \quad \sum_{m=1}^{\infty}\left|\left(x_m, x\right)\right|=\infty \ \ \forall  x\in H, x\neq 0$$
Then the set of all convergent series
$
\sum_{m=1}^{\infty} a_m x_m, \quad\left|a_m\right|=1, m=1,2, \ldots
$
is dense in $H$.
\end{graybox}
  The proof of this result is a testimony of the power of the Hahn-Banach theorem\endnote{Besides the Hahn-Banach theorem one uses the following Lemma (see \cite{Laurincikas} Lemma 1.15) : Let $x_1, \ldots, x_n \in H$ and let $b_1, \ldots, b_n$ be complex numbers with $\left|b_j\right| \leqslant 1, j=1, \ldots, n$. Then there exist complex numbers $a_1, \ldots, a_n$ with $\left|a_j\right|=1, j=1, \ldots, n$, such that
$$
\left\|\sum_{j=1}^n a_j x_j-\sum_{j=1}^n b_j x_j\right\|^2 \leqslant 4 \sum_{j=1}^n\left\|x_j\right\|^2
$$}, using the convex set of sums $\sum b_m x_m$ where $\vert b_m\vert \leq 1$.\newline      
    Voronin's theorem has been extended to joint universality for Dirichlet $L$-functions and other classes of $L$-functions, revealing this as a general phenomenon for functions with Euler products. The essence of the proof is to show that, as elements of a suitable Bergman/Hardy-type space $H$, the system $\left\{\log \left(1-p^{-s}\right)\right\}$ where $p$ runs through the primes  satisfies both conditions of the above theorem. This allows one to approximate the logarithm of the function $f(s)$ (which makes sense since $K$ is simply connected).

\section{ A Century and a Half of Theory Building towards RH}
Since Riemann's 1859 memoir ``Über die Anzahl der Primzahlen unter einer gegebenen Grösse,'' the hypothesis he formulated concerning the zeros of the zeta function has inspired some of the most profound developments in mathematics. 
The pursuit of the Riemann Hypothesis (RH), whether explicitly or as a motivating ideal, has led to the emergence of entire theories and reshaped multiple domains. This section provides an overview of major mathematical frameworks and ideas developed in this ongoing quest.

\subsection{Harmonic Analysis and Functional Analysis}

\subsubsection{Hilbert spaces and spectral theory:} The proof of Voronin's  theorem discussed above shows the power of Hilbert space techniques. Hilbert's inclusion of RH as his 8th problem (1900) coincided with the development of spectral theory. The Hilbert-Pólya conjecture (1910s) suggests the existence of  a self-adjoint operator $H$ such that the non-trivial solutions of  $\zeta(1/2 + it) = 0$ are the eigenvalues of $H$. This would immediately imply RH since eigenvalues of self-adjoint operators are real.

    \subsubsection{Scattering theory and spectral interpretation:} The scattering approach, developed by Faddeev and Pavlov, interprets the zeros of $\zeta(s)$ as resonances (poles of the scattering matrix) of a certain quantum mechanical system. The Berry-Keating conjecture \cite{BerryKeatingChaos} suggests connections to classical chaotic Hamiltonian systems whose quantization might yield the desired operator. Specific proposals include $H =\frac 12( xp+px)$ (where $x$ is position and $p$ is momentum) with appropriate boundary conditions, though rigorous constructions remain elusive. The connection to quantum chaos is reinforced by the statistical properties of the zeros matching those of quantum chaotic systems.

\subsection{Algebraic and Arithmetic Geometry}

As every mathematician knows, one of the time-honored strategies when faced with a difficult problem is to enlarge its scope: to generalize the question and then examine special cases where essential features can be more clearly discerned. Within this setting, Andr\'e Weil achieved a decisive breakthrough by resolving the analogue of the hypothesis for global fields of finite characteristic. 

\subsubsection{Weil's proof for function fields:} Weil's 1940s proof of RH for curves over finite fields used intersection theory on the product $C \times C$ of a curve with itself. For a smooth projective curve $C$ over $\mathbb{F}_q$, the zeta function
    \[
    Z_C(t) = \exp\left(\sum_{n=1}^{\infty} \frac{N_n}{n}t^n\right)
    \]
    where $N_n = |C(\mathbb{F}_{q^n})|$, satisfies a functional equation and has the form $Z_C(t) = P(t)/((1-t)(1-qt))$ where $P(t)$ is a polynomial whose roots have absolute value $q^{-1/2}$. This is the analogue of RH. Weil's proof used the intersection pairing on divisors and the Hodge index theorem, establishing a geometric approach to RH-type problems.  The successive simplifications of Weil's proof  by Mattuck-Tate and Grothendieck represent a paradigmatic example of how mathematical understanding evolves through the development of more powerful theoretical frameworks. Weil's original achievement required  technical virtuosity to overcome the absence of key theoretical tools. The subsequent simplifications demonstrated how the availability of general results—the Riemann-Roch theorem for surfaces\endnote{It is of the form    
    \begin{equation}
    \sum_0^2(-1)^j {\rm dim}\, H^j(X,\mathcal O(D))=\frac 12 D.(D-K)+\chi(X)
\end{equation}
where $D$ is a divisor, $\mathcal O(D)$ its sheaf of sections, $\chi(X)$ is the arithmetic genus}
 and systematic intersection theory—could transform specialized constructions into applications of standard techniques. While the classical Riemann-Roch theorem for surfaces originated with Castelnuovo in the 1890s, the sheaf-theoretic version that enabled the Mattuck-Tate simplification was developed by Serre, and the intersection theory suitable for this setting represents a collaborative effort involving André Weil's foundational 1946 work, Claude Chevalley's local theory, Serre's 1958 algebraic formulation, and subsequent developments that created the modern framework of intersection multiplicities using Tor functors and sheaf cohomology.

   \subsubsection{Grothendieck's schemes and étale cohomology:} Grothendieck revolutionized algebraic geometry by introducing schemes and  étale cohomology \cite{ArtinGrothendieckVerdier}, partly motivated by the goal of proving the Weil conjectures in arbitrary dimension.  The étale cohomology groups $H^i_{\text{ét}}(X, \mathbb{Q}_{\ell})$ for a variety $X$ over $\mathbb{F}_q$ carry an action of the Frobenius morphism, and the zeta function is expressed as
    \[
    Z_X(t) = \prod_{i=0}^{2\dim X} P_i(t)^{(-1)^{i+1}}
    \]
    where $P_i(t) = \det(1 - tF | H^i_{\text{ét}}(X, \mathbb{Q}_{\ell}))$. Deligne's proof (1974) of the Weil conjectures showed that the eigenvalues of Frobenius on $H^i$ have absolute value $q^{i/2}$, establishing RH for varieties over finite fields.\newline
      
   What is particularly inspiring in the study of the \'etale site of $\mathrm{Spec}(\mathbb{Z})$ is the striking analogy, due to David Mumford and Barry Mazur, between the role of prime ideals in arithmetic geometry and the role of knots in the three-dimensional sphere. In this perspective, each prime can be viewed as the analogue of a knot, and the intricate way in which primes are intertwined within $\mathrm{Spec}(\mathbb{Z})$ endowed with the \'etale topology mirrors the way knots can be interlaced in $S^3$. This vision opened a new topological imagination for number theory, suggesting that the geometry of primes might be understood through concepts closer to low-dimensional topology than to classical analysis. At the heart of this analogy lies Grothendieck’s far-reaching extension of Galois theory from fields to general schemes, accomplished through the introduction of the \'etale fundamental group. It provides the natural language to interpret the arithmetic of primes in terms of covering spaces, thereby giving a conceptual framework in which Mazur’s and Mumford’s analogy could take shape and develop~\cite{Mazur1973}. We shall come back to this point briefly in \S \ref{adelic}.

   \subsubsection{Motives} Grothendieck's vision of motives aims to provide a universal cohomology theory that explains the common patterns in various cohomology theories (de Rham, étale, crystalline, etc.). In this framework, zeta and $L$-functions are associated to motives, and their analytic properties (including potential RH) would follow from geometric properties of the motives. The motivic formalism suggests that there should be a cohomological interpretation of the classical RH, viewing $\mathrm{Spec}(\mathbb{Z})$ as an arithmetic curve and seeking an appropriate cohomology theory. While the full theory of mixed motives remains conjectural, special cases like Artin motives and motives of modular forms have been extensively developed.
\subsection{Automorphic Forms and Representation Theory}

In Andr\'e Weil’s work, the Riemann Hypothesis is naturally recast within the broader framework of global fields. He proved the analogue of the hypothesis for global fields of finite characteristic, and his book \emph{Basic Number Theory} (Springer, 1967) illustrates how powerful it is to approach number theory from this unifying perspective. A fundamental insight is that, although global fields themselves are countable and discrete, each admits a natural companion: a locally compact, non-discrete ring that contains the global field as a discrete co-compact subgroup. This companion is the ring of ad\`eles, first introduced in the context of class field theory. The ad\`ele ring provides a rich harmonic-analytic structure, and Weil’s reformulation of the explicit formulas in this setting represents a decisive step toward harnessing that structure as a tool in analytic number theory. Weil’s book \emph{Basic Number Theory}  was the first systematic exposition of the use of ad\`eles and id\`eles in number theory. It not only unified class field theory but also laid the foundation for later developments in automorphic forms and representation theory to which we now turn.

\subsubsection{Langlands Program:} This extension of class field theory to the non-abelian case is a vast web of conjectures, initiated by Langlands in the 1960s, it predicts deep connections between Galois representations and automorphic representations. For $L$-functions, it predicts that every motivic $L$-function (coming from algebraic geometry) equals an automorphic $L$-function (coming from harmonic analysis on adelic groups). The proven functional equation is actually one of the key achievements of the automorphic theory, establishing that automorphic $L$-functions have the right analytic structure to even formulate the Generalized Riemann Hypothesis. Since automorphic $L$-functions are expected to satisfy  RH, this would imply RH for all motivic $L$-functions. The program has achieved spectacular successes (Wiles' proof of Fermat's Last Theorem used the Taniyama-Shimura conjecture, now a theorem) but the general conjectures remain open.

   \subsubsection{Modular forms and L-functions:} The theory of modular forms provides concrete examples of $L$-functions with analytic continuation and functional equations. For a weight $k$ cusp form $f(z) = \sum_{n=1}^{\infty} a_n e^{2\pi i nz}$ on $\Gamma_0(N)$, the associated $L$-function
    \[
    L(s,f) = \sum_{n=1}^{\infty} \frac{a_n}{n^s}
    \]
    satisfies a functional equation relating $s$ to $k-s$. The Ramanujan-Petersson conjecture (proved by Deligne) bounds $|a_p| \leq 2p^{(k-1)/2}$, which is the analogue of RH at the level of Euler factors. These examples clarify the expected shape of general $L$-functions and their connection to geometric objects.

   \subsubsection{Selberg trace formula and zeta functions:} Selberg's trace formula \cite{SelbergTrace1,Hejhal} for discrete groups $\Gamma$ acting on hyperbolic spaces relates spectral data (eigenvalues of the Laplacian) to geometric data (lengths of closed geodesics). For $\Gamma = \mathrm{PSL}(2,\mathbb{Z})$, this yields the Selberg zeta function
    \[
    Z_{\Gamma}(s) = \prod_{p \text{ primitive}} \prod_{k=0}^{\infty} (1 - e^{-(s+k)l(p)})
    \]
    where the product is over primitive closed geodesics $p$ with length $l(p)$. This function has zeros at the eigenvalues of the Laplacian, providing a complete spectral interpretation. The analogy with the Riemann zeta function is striking: both have Euler products, functional equations, and their zeros encode important information (eigenvalues vs. prime distribution). The Selberg  trace formula for Riemann surfaces of finite area, acquires additional terms which make it look \eg in the case of $X=H/PSL(2,\Z)$ (where $H$ is the upper half plane with the Poincar\'e metric) even more similar to the explicit formulas, since the parabolic terms now involve explicitly the sum $$2\sum_{n=1}^\infty \frac{\Lambda(n)}{n}g(2\log n)$$  
(for a test function $g$) to be compared with   the $\Lambda(n)$ terms in the explicit formulas 
$$
-2\sum_{n=1}^\infty \frac{\Lambda(n)}{n^{\frac 12}}g(\log n)
$$
There is however a striking difference which is that these terms occur with a positive sign instead  of the   negative sign, as discussed in \cite{Hejhal} \S 12. This discussion of the minus sign was extended to the case of the semiclassical limit of Hamiltonian
systems in physics in \cite{BerryKeatingChaos}.

\subsection{Random Matrix Theory and Quantum Chaos}

\subsubsection{Montgomery's pair correlation:}     
    In 1973, Hugh Montgomery \cite{Montgomery,MontgomeryPair}  conjectured a striking statistical property of the nontrivial zeros of the Riemann zeta function on the critical line. Letting $\gamma$ denote ordinates of zeros $\rho = \tfrac{1}{2} + i\gamma$ with $0 < \gamma < T$, he studied the two-point correlation function of the rescaled spacings between zeros. Specifically, he conjectured that for $0 < a < b$, and $N(T)=\sum_{0<\gamma \leq T} 1$,
$$
\lim_{T \to \infty} \frac{1}{N(T)} \# \left\{ (\gamma, \gamma’) : 0 < \gamma, \gamma’ < T,\; \frac{2\pi a}{\log T} \leq |\gamma - \gamma’| \leq \frac{2\pi b}{\log T} \right\} = \int_a^b \left( 1 - \left( \frac{\sin \pi u}{\pi u} \right)^2 \right) du,
$$
 This conjectured density reflects a repulsion between zeros and suggests they are not randomly distributed but exhibit a regularity akin to that of eigenvalues of large random Hermitian matrices. In a widely recounted anecdote, Montgomery mentioned his formula during a visit to the Institute for Advanced Study, where the physicist Freeman Dyson, upon seeing it, immediately recognized it as the pair correlation function for eigenvalues in the Gaussian Unitary Ensemble (GUE) of random matrix theory. This was verified numerically by Odlyzko to extraordinary precision. The connection suggests that the zeros behave like eigenvalues of a random Hermitian matrix, pointing to an underlying chaotic quantum system.
 This serendipitous encounter laid the groundwork for a deep and fruitful connection between number theory and quantum statistical mechanics.
    
    \subsubsection{Odlyzko's statistics:} Andrew Odlyzko's computations in the 1980s-1990s revolutionized our understanding of the statistical properties of zeta zeros. Computing millions of zeros at very large heights (around the $10^{20}$-th zero), he found:
    \begin{itemize}
        \item Nearest-neighbor spacings follow the GUE distribution with extraordinary accuracy
        \item Higher correlation functions also match random matrix predictions
        \item Local statistics are universal, but global statistics show number-theoretic fluctuations
    \end{itemize}
    His data provided compelling evidence for the GUE hypothesis and inspired much of the subsequent work connecting number theory to random matrix theory. The agreement extends to fine details like the variance of the number of zeros in intervals and the distribution of the arguments of $\zeta(1/2 + it)$.\newline
    
    To compare the statistical properties of Riemann zeta zeros with eigenvalues of random matrices from the Gaussian Unitary Ensemble (GUE), one must address a fundamental incompatibility in their densities. The eigenvalues of an $N \times N$ GUE matrix have constant mean spacing $\pi/\sqrt{N}$ along the real axis (after appropriate normalization). In contrast, the zeros $\rho_n = 1/2 + i\gamma_n$ of the Riemann zeta function have imaginary parts $\gamma_n$ whose density increases logarithmically: the number of zeros with $0 < \gamma < T$ is approximately $\frac{T}{2\pi}\log\frac{T}{2\pi}$, giving a local mean spacing near height $T$ of approximately $\frac{2\pi}{\log(T/2\pi)}$. 

To make a meaningful comparison, one must rescale the zeros locally. Specifically, when studying zeros near height $T$, one considers the rescaled spacings
$$\tilde{\gamma}_n = \frac{\log(T/2\pi)}{2\pi}(\gamma_{n+1} - \gamma_n)$$
which have mean value 1. One then compares these rescaled spacings with the eigenvalue spacings of GUE matrices with $N \sim \log\frac{T}{2\pi}$. After this height-dependent rescaling, one finds remarkable agreement in local statistics such as nearest-neighbor spacing distributions and $n$-point correlation functions. However, this rescaling procedure highlights that any correspondence between zeta zeros and random matrix eigenvalues must be inherently local and cannot arise from a simple, fixed spectral operator.

 \subsubsection{Quantum chaos:} The field of quantum chaos studies the quantum mechanics of classically chaotic systems. The Berry-Tabor conjecture states that quantum systems with integrable classical limits have Poisson statistics for their eigenvalues, while those with chaotic classical limits follow random matrix statistics (GUE for time-reversal invariant systems). The appearance of GUE statistics for zeta zeros suggests they arise from quantizing a chaotic classical system \cite{BerryKeatingChaos}. Specific proposals include billiards in certain domains or more abstract dynamical systems.

   \subsubsection{Katz-Sarnak theory}\label{sectKS} Building upon Montgomery's discovery  and Odlyzko's 
computational confirmation  of the connection 
between Riemann zeta zeros and random matrix statistics, Nicholas Katz and 
Peter Sarnak extended this correspondence to entire families of $L$-functions 
and established a systematic theoretical framework for understanding their 
statistical behavior \cite{KatzSarnakBook,KatzSarnakBull}. Their 
groundbreaking insight was that different families of $L$-functions exhibit 
universal statistical properties governed by one of the classical matrix 
groups---unitary $U(N)$, orthogonal $O(N)$, or symplectic $Sp(2N)$---with the 
specific symmetry type determined by the arithmetic and geometric structure of 
the family itself \cite{KatzSarnakBook,Duenez2004}. The Katz--Sarnak density 
conjecture predicts that, as the analytic conductors of $L$-functions in a 
family tend to infinity, the distribution of their normalized low-lying zeros 
near the critical point $s = 1/2$ converges to the scaling limits of 
eigenvalues clustered near $1$ in the corresponding random matrix ensemble 
\cite{KatzSarnakBook,KatzSarnakBull,KeatingSnaith2000}. Most remarkably, 
they proved their conjectures rigorously for families of $L$-functions over 
finite fields by connecting the zeros to Frobenius eigenvalues and applying 
Deligne's equidistribution theorem \cite{KatzSarnakBook,GordonSarnak1990,
DeligneWeilI,DeligneWeilII}, thereby providing the first complete verification 
of random matrix universality in a number-theoretic context. This work not only 
unified previously disparate phenomena under a single conceptual framework but 
also provided powerful new tools for studying classical problems in analytic 
number theory, from the distribution of primes to the arithmetic of elliptic 
curves.

\subsubsection{Moments and unitary matrices (Keating–Snaith):}
The asymptotic behavior of moments of the zeta function is conjectured (see \cite{Conrey84}) to be of the form :
\[
\frac 1T\int_0^T |\zeta(1/2 + it)|^{2k} dt \sim c_k (\log T)^{k^2}
\]
for some constant $c_k$.  Hardy and Littlewood proved this when $k=1$ in 1918 and Ingham proved it in 1926 \cite{Ingham2} for $k=2$, one has $c_1=1$, and $c_2=\left(2 \pi^2\right)^{-1}$.  These are the only values of $k$ where this conjecture is proven, and there was for a long time not even a conjectural value of $c_k$ for any other value of $k$. Conrey and Ghosh \cite{Conrey3} conjectured in 1992 the answer for $k=3$, and Conrey and Gonek \cite{Conrey4} conjectured the answer for $k=4$ in 1998. These two conjectures are of the form $c_k=a_k g_k/(k^2!)$ where as explained below $a_k$ is a well understood arithmetic factor, but the constants $g_3=42$ and $g_4=24024$ remained mysterious. In the 1998 Vienna conference on RH,  Keating and Snaith \cite{KeatingSnaith2000,KeatingSnaithMoments} announced their breakthrough of using random matrix theory to guess the value of $g_k$. The new insight of Keating and Snaith was recognizing that the constant $f_k=g_k/(k^2!)$ naturally arises from 
random matrix theory. The general conjecture then takes the form $c_k=a_kf_k$ where
\begin{itemize}
\item $f_k$ comes from random matrix theory (the $k$-th moment of $|\det(I - U)|^2$ for $U$ a random unitary matrix)
\item $a_k$ is an arithmetic factor which captures the arithmetic complexity through an Euler product, $$a_k=\prod_p\left(1-\frac{1}{p}\right)^{k^2} \sum_{m=0}^{\infty}\left(\frac{\Gamma(m+k)}{m!\Gamma(k)}\right)^2 p^{-m}$$
\end{itemize}
 The random matrix part $f_k$ can be expressed as:
\[
f_k = \prod_{j=0}^{k-1} \frac{j!}{(j+k)!}\implies g_k=(k^2)!\,f_k
\]
 This separation reveals
  a beautiful principle: the ``universal'' behavior comes from random matrix theory, while the arithmetic specifics of the zeta function are encoded separately in $a_k$. This opened a new paradigm of using random matrix models to predict number-theoretic results, later extended to other $L$-functions and their derivatives \cite{ConreySnaithRatios,GonekHughesKeatingHybrid,Snaith}.

\subsection{Noncommutative Geometry}

\subsubsection{Connes' trace formula} Connes developed a trace formula in the context of noncommutative geometry\footnote{The origin came from  a system of quantum statistical mechanics exhibiting phase transitions and having the Riemann zeta function as partition function \cite{BC}} that recovers Weil's explicit formula for the zeros of $\zeta(s)$ and of $L$-functions \cite{Co-zeta, CMbook}. The key insight is to work with the noncommutative space of adele classes $\mathbb{A}/\mathbb{Q}^*$, where $\mathbb{A}$ is the ring of adeles. On this space, there is a natural action of the idele class group $\operatorname{GL}(\mathbb{A})/\mathbb{Q}^*\simeq \hat\Z^*\times \R_+^*$. The trace formula relates:
    \begin{itemize}
        \item Spectral side: contributions from the zeros of $L$-functions with Grössencharakter.
        \item Geometric side: contributions from places of $\Q$.
    \end{itemize}
    In the  beginning the trace formula was stated for general global fields as a conjecture equivalent to the validity of the analogue of RH for all $L$-functions with Grössencharakter, but the semilocal form of the trace formula was already proved in \cite{Co-zeta}, and it admits as simple corollaries the global versions allowing for the presence of non-critical zeros \cite{Meyer} (see also \cite{CMbook}). 
    This provides a spectral realization where the zeros appear as an absorption spectrum which accounts for the minus sign in the explicit formulas as compared to the Selberg trace formula. \newline
    These approaches aim to geometrize the RH using tools from operator algebras and noncommutative differential geometry. 
   An ergodic quotient such as the adele class space $\A_K/K^\times$ of a number field $K$ is encoded by the cross product of the algebra of functions on $\A_K$ by the action of $K^\times$. The spectral realizations mentioned above involve the Hochschild homology of the cross product of the Bruhat-Schwartz algebra, while the topology of $\A_K/K^\times$ is encoded by the $C^*$-algebra cross product. A recent result \cite{Bruce} shows that this topological encoding is faithful, non-isomorphic number fields give non-isomorphic $C^*$-algebras.
    \subsubsection{Knots, primes and class field theory}\label{adelic}
    For the global field $K=\Q$, the projection  of the adele class space
    $$
    \pi:Y_\Q=\mathbb{Q}^*\backslash \mathbb{A}\to \mathbb{Q}^*\backslash \mathbb{A}/\hat\Z^*=X_\Q
    $$
    on the sector $X_\Q$ corresponding to $\zeta$ displays the class field theory counterpart of the Mazur-Mumford analogy between knots and primes mentioned above. To each prime $p$ corresponds a periodic orbit $C_p\subset X_\Q$ of length $\log p$. 
    Let then ${\rm F r o b}_{p}\in \pi_1^{e t}(\Spec(\F_p))$ be the canonical generator of the etale fundamental group and $\Z_{(p)}$ the ring $\Z$ localized at $p$. One then has \cite{ccknots}, 
     \begin{enumerate}
         \item[(i)] The inverse image $\pi^{-1}(C_p)\subset Y_\Q$ of the periodic orbit $C_p$  is canonically isomorphic to the mapping torus of the multiplication by $r^*\left\{{\rm F r o b}_{p}\right\}$  in the abelianized \'etale fundamental group $\pi_1^{e t}(\Spec \, \Z_{(p)})^{ab}$.
          \item[(ii)] The canonical isomorphism in $(\rm i)$ is equivariant for the action of the idele class group.
         \item[(iii)] The monodromy of the periodic $C_p$ in $\pi^{-1}(C_p)\subset Y_\Q$ is equal to the natural map $$r^*:\pi_1^{e t}\left(\operatorname{Spec} \mathbb{F}_p\right) \rightarrow \pi_1^{e t}(\operatorname{Spec} \mathbb{Z}_{(p)})^{a b}$$ and determines the linking of the prime $p$ with all other primes.
         \end{enumerate}
         This result shows that the adele class space plays the role of a class field theory counterpart of the abelian etale covers involved in Grothendieck theory. 
         Quite remarkably the space $X_\Q$ admits a topos theoretic incarnation as the scaling site naturally endowed with a structure sheaf of characteristic one (see \cite{Essay} and the references there).

\subsection{p-adic and Motivic L-functions}

\subsubsection{p-adic L-functions and Iwasawa theory:} Kubota and Leopoldt constructed $p$-adic analogues of the Riemann zeta function, which are $p$-adic analytic functions that interpolate special values of the classical zeta function. For each prime $p$, the $p$-adic zeta function $\zeta_p(s)$ is a $p$-adic analytic function on $\mathbb{Z}_p$ satisfying
    \[
    \zeta_p(1-k) = (1-p^{k-1})\zeta(1-k)
    \]
    for positive integers $k$. Iwasawa theory studies these functions in towers of cyclotomic extensions, relating them to ideal class groups and units. The Main Conjecture of Iwasawa theory (proved by Mazur-Wiles) connects $p$-adic $L$-functions to Selmer groups of Galois representations. While there is no direct $p$-adic analogue of RH, the growth and zeros of $p$-adic $L$-functions are intimately connected to deep arithmetic phenomena.

  \subsubsection{Motivic L-functions and Bloch–Kato conjectures:} These far-reaching conjectures relate special values of $L$-functions to arithmetic invariants. For a motive $M$ with $L$-function $L(M,s)$, the conjectures predict that:
    \begin{itemize}
        \item The order of vanishing of $L(M,s)$ at integer points equals the rank of certain $K$-groups or Selmer groups
        \item The leading coefficient is related to regulators, periods, and Tamagawa numbers
    \end{itemize}
    Special cases include the Birch and Swinnerton-Dyer conjecture (for elliptic curves) and Beilinson's conjectures. These place RH in a broader arithmetic context where zeros of $L$-functions encode geometric and arithmetic information. The conjectures suggest that RH is part of a vast network of relationships between analysis, algebra, and geometry.

\subsection{Computational and Experimental Mathematics}

\subsubsection{High-precision computations} The numerical verification of RH has a rich history:
    \begin{itemize}
        \item Riemann (1859): computed the first few zeros
        \item Gram (1903): 15 zeros
        \item Backlund (1914): 79 zeros
        \item Hutchinson (1925): 138 zeros
        \item Titchmarsh (1935-1936): 1,041 zeros using the Riemann-Siegel formula
        \item Turing (1950):  verified 1,104 zeros using the Manchester Mark 1
        \item Lehmer (1956): 25,000 zeros using electronic computers
        \item Rosser, Yohe, Schoenfeld (1968): 3,500,000 zeros
        \item van de Lune, te Riele, Winter (1986): 1,500,000,000 zeros
        \item Gourdon and Demichel (2004): first $10^{13}$ zeros
        \item Platt (2021): verification up to height $3 × 10^{12}$, \cite{Platt}.
    \end{itemize}
    Modern computations use the Riemann-Siegel formula with sophisticated error bounds and the Turing method for rigorous verification. The 1988 publication of the Odlyzko-Schönhage algorithm was a major methodological breakthrough. These massive computations provide overwhelming evidence for RH while also testing for exceptional phenomena \cite{Edwards,IvicBook}.

\section{Equivalent Formulations}

While the previous sections surveyed theories built in pursuit of the Riemann Hypothesis, there exists a parallel tradition of discovering equivalent formulations of RH itself. These reformulations, ranging from elementary number theory to functional analysis, reveal the hypothesis's deep connections across mathematics and sometimes suggest new angles of attack.  There are so many of these equivalent formulations that there is a whole website which is designed to account for them. The main significance that one can find in all these formulations is that the elementary formulations, for instance, the Robbins criterion or the Lagarias  criterion, show that at the logical level RH has a very particular logical  status which will be discussed in details below.

\subsection{Weil's Positivity Criterion}\label{sectweilpos}
 RH is equivalent to the positivity of certain distributions constructed from the zeros, connecting to Weil's approach via the explicit formula.
The difficulty of solving the Riemann Hypothesis in the analytic formulation, is often  
primarily attributed to the infinite number of terms within the Euler product.
\begin{equation}\label{eulerproduct1}
 	\zeta(s) = \prod_p (1- p^{-s})^{-1}
 \end{equation} 
However, contrary to this widespread belief,  there exists a property 
$P(n)$, involving only
the Euler factors for primes smaller than $n$, and
 whose validity for all $n$ is equivalent to RH. 
This is derived  from  Weil's positivity criterion which involves the quadratic
 form $QW$ defined using the Riemann-Weil explicit formulas applied to test functions with support in a compact symmetric interval. After a minor change of notations the Mellin transform becomes the Fourier transform for the group $\R_+^*$ whose Pontjagin dual is identified with the additive group $\R$, and the explicit formula takes the following  form, similar to \eqref{bombieriexplicit1}, \eqref{bombieriexplicit2} but with notations adapted to the involutive convolution algebra of the group $\R_+^*$,
 $$\begin{aligned} & \widehat{f}\left(\frac{i}{2}\right)-\sum_{\frac{1}{2}+i s \in Z} \widehat{f}(s)+\widehat{f}\left(-\frac{i}{2}\right)=\sum_v W_v(f) \\ & \widehat{f}(s):=\int_0^{\infty} f(x) x^{-i s} d^* x, \quad d^* x=\frac{d x}{x}\end{aligned}$$
where the local contributions are now given by 
\begin{equation}\label{wp}W_p(f):=(\log p) \sum_{m=1}^{\infty} p^{-m / 2}\left(f\left(p^m\right)+f\left(p^{-m}\right)\right)
 \end{equation}
 and for the archimedean place
\begin{equation}\label{arch}
\begin{aligned} W_{\mathbb{R}}(f):= & (\log 4 \pi+\gamma) f(1) \\ & +\int_1^{\infty}\left(f(x)+f\left(x^{-1}\right)-2 x^{-1 / 2} f(1)\right) \frac{x^{1 / 2}}{x-x^{-1}} d^* x\end{aligned}	
\end{equation}
The key result of André Weil is the equivalence
$$R H \Longleftrightarrow \sum_v W_v\left(g * g^*\right) \leq 0,\ \forall g, \  \widehat{g}\left( \pm \frac{i}{2}\right)=0$$
where $g \in C_c^{\infty}\left(\mathbb{R}_{+}^*\right)$ is smooth with compact support and  $g^*(x):=\bar{g}\left(x^{-1}\right)$.
The key point of this equivalence is that the sum on the right-hand side, when evaluated on a test function $g$ with compact support, involves only finitely many primes (since $W_p$ vanishes on functions with support in $\left(p^{-1}, p\right)$ ). Thus, even though the Riemann Hypothesis is concerned with the asymptotic distribution of the primes, the equivalent formulation  only involves finitely many primes at a time. In \cite{yoshida} H. Yoshida proved the following result (Theorem 1 in that paper)
\begin{graybox}
For any smooth, positive definite function $f$ with support in the interval $(1 / 2,2)$ and whose Fourier transform vanishes at $\pm \frac{i}{2}$ one has: $W_{\infty}(f) \geq 0$ where $W_{\infty}:=-W_{\mathbb{R}}$.
\end{graybox}
The proof is a numerical analysis of the positivity of the Weil functional $W_{\infty}$ restricted to the interval $\left(\frac{1}{2}, 2\right)$, and therefore it does not provide any conceptual reason for this positivity that would have a chance to continue to hold when primes are involved.

 \subsection{Beurling–Nyman criterion:} This remarkable reformulation (Beurling 1955, Nyman 1950) states that RH is equivalent to the density  in $L^2(0,1)$ of  the linear combinations $\sum_1^n c_\nu\rho_{\theta_\nu}$, $0<\theta_\nu \leq 1$, $\sum_1^n c_\nu \theta_\nu=0$ where
    \[
    \rho_{\theta}(x) =\{ \frac{\theta}{x}\} \quad \text{for } x \in (0, 1)
    \]
    and $\{y\}$ denotes the fractional part of $y$. This transforms RH into a completeness problem in functional analysis. Later work by Báez-Duarte and others has provided quantitative versions, showing that the rate of approximation in this criterion is related to the distribution of zeros off the critical line.
    
    \subsection{Li's criterion}
   Xian-Jin Li (1997) proved \cite{XJL} that RH is equivalent to the positivity of the numbers
$$\lambda_n = \sum_{\rho} \left(1 - \left(1 - \frac{1}{\rho}\right)^n\right)$$
for all $n \geq 1$, where the sum runs over all non-trivial zeros $\rho$ of $\zeta(s)$.

\subsection{Elementary Number-Theoretic Formulations}

In 1913 Grönwall proved in his paper "Some asymptotic expressions in the theory of numbers" that $$\limsup _{n \rightarrow \infty} \frac{\sigma(n)}{n \log \log n}=e^\gamma, \  \ \gamma=\lim _{n \rightarrow \infty}\left(\sum_{k=1}^n \frac{1}{k}-\log n\right)$$ where $\sigma(n)$ is the sum of divisors function and $\gamma$ is the Euler-Mascheroni constant. Building on Gronwall's foundation, Srinivasa Ramanujan made a crucial connection between the  function $\sigma(n)$ and the Riemann Hypothesis in his work on highly composite numbers. Ramanujan proved that if the Riemann Hypothesis is true, then the inequality $\sigma(n) / n<e^\gamma \log \log n$ holds for all sufficiently large positive integers $n$. This represented the first direct link between the truth of the Riemann Hypothesis and bounds on arithmetic functions.

\subsubsection{Robin's Criterion (1984)} Guy Robin proved that RH is equivalent to the inequality
$$\sigma(n) < e^\gamma n \log \log n$$
for all $n > 5040$, where $\sigma(n)$ is, as above, the sum of divisors function and $\gamma$ is the Euler-Mascheroni constant. This remarkable result translates the complex-analytic RH into purely arithmetic terms.

\subsubsection{Lagarias' Criterion (2002)} Jeffrey Lagarias refined Robin's criterion in the form
$$RH\iff \sigma(n) < H_n + e^{H_n} \log H_n, \  \ \forall n \geq 1$$
 where $H_n = 1 + 1/2 + \cdots + 1/n$ is the $n$-th harmonic number. This formulation has the aesthetic advantage of holding for all positive integers.\newline
 The reformulation of the Riemann Hypothesis as a universal statement over decidable arithmetic properties places it precisely within the class of statements that Hilbert was hoping to consider as "provable" if true and that Gödel's incompleteness theorems identify as potentially "true but unprovable" ! Chaitin's algorithmic information theory results demonstrate that such unprovability becomes increasingly prevalent as statement complexity grows. This connection illuminates fundamental limitations in formal mathematical systems and provides a striking example of how even central problems in number theory intersect with the deepest questions about the nature of mathematical truth and proof.
   Chaitin's work reveals that mathematical truth and provability operate in fundamentally different domains from an information-theoretic perspective. Most mathematical truths have high descriptive complexity and contain more information than can be extracted from finite axiom systems through deductive processes. This creates a vast landscape of true statements that remain forever beyond the reach of formal proof, with provable statements representing an infinitesimal fraction of mathematical reality as complexity increases. 
Estimating the algorithmic complexity of the Lagarias criterion reveals fundamental computational barriers in verifying the Riemann Hypothesis through direct arithmetic checking. The dominant complexity arises from integer factorization requirements for computing the sum of divisors function, placing individual verification instances outside polynomial-time computation for general integers. Harmonic number computation and transcendental function evaluation contribute additional complexities that compound the overall computational costs.

These complexity limitations highlight both the elegance and computational intractability of the Lagarias criterion as a verification approach for the Riemann Hypothesis. While the criterion successfully transforms a deep analytic statement into elementary arithmetic, the computational costs of verification reflect the profound difficulty of the underlying mathematical problem. The analysis demonstrates how computational complexity theory provides essential insights into the practical feasibility of mathematical verification approaches and the fundamental relationship between mathematical truth and computational resources.

\newpage
\section{ A Letter to Professor Bernhard Riemann}
 I was invited to give a talk in Varese, on June 4th of 2025, in the Villa Toeplitz, where  the Riemann International School of Mathematics is located. In order to prepare this talk, I did, two days before, a pilgrimage in Selasca,  the place where Riemann died on  July 20th of 1866. 

Given the place where I had to give the talk, I took the following challenge:

\medskip
\centerline{\it What could I tell to Riemann that would surprise him} \centerline{\it and that would give him confidence that his hypothesis is true?}
\medskip
 So what I am going to write now is a letter to him. I will address to him as "Master". The requisite of course is that it is forbidden to use mathematical  notions unknown to him or  not  easily understandable by him. In fact, I will restrict myself  to those which he  used in his work. After a  preliminary remark, the really new stuff is in the second part of the letter. 
 \bigskip
 %\subsection{The Letter}
 \begingroup
%\fontfamily{EBGaramond-LF}\selectfont
\newenvironment{ancientletter}{%
  \fontfamily{ppl}\selectfont Palatino
  \setlength{\parindent}{1em}
  \setlength{\parskip}{0.5em}
  \large % slightly larger text
}{%
  \normalfont\normalsize
}

 \begin{center}
{\Large\scshape To Professor Bernhard Riemann}
\end{center}
 
 \begin{large}%\centerline{\bf Master}
 
  \bigskip
 
 The beautiful formula that you proved in your paper, namely the formula which allows one to compute the number of prime numbers less than $x$ in terms of the zeros of the zeta function, has been unfortunately badly written in textbooks thus displaying a profound misunderstanding of the middle term of the formula.  You let 
  $\pi'(x)$ be the number of primes strictly less than $x$ with $\frac{1}{2}$ added when $x$ is prime and  found for the counting function 
$$
f(x):=\sum \frac{1}{n} \pi'\left(x^{\frac{1}{n}}\right),
$$
the following formula involving the integral logarithm function $\operatorname{Li}(x)=\int_0^x \frac{d t}{\log t}$,
\begin{equation}\label{piriemann}
f(x)=\operatorname{Li}(x)-\sum_\alpha\left(\operatorname{Li}\left(x^{\frac{1}{2}+\alpha i}\right)+\operatorname{Li}\left(x^{\frac{1}{2}-\alpha i}\right)\right)+\int_x^{\infty} \frac{1}{t^2-1} \frac{d t}{t \log t}-\log 2
 \end{equation}
 You were very careful at the end of your paper in dealing with multivalued functions, but the mathematical language has evolved and, nowadays, one only uses functions with  a definite univalent meaning, while the use of multivalued functions is no longer commonly admitted.\newline 
  In standard textbooks, for instance the classic book of Edwards\endnote{It is not clear if there is a hidden convention in this book, but the following formula written in section 1.15 "The term involving the roots $\rho$" 
  $$
  \int_{C^+}\frac{t^{\beta-1}}{\log t}dt=\int_0^{x^\beta}\frac{du}{\log u}
  $$
  (with the comment "where the second integral is over a path which passes above the singularity at $u=1$")
  equates two terms which cannot be equal since the integral on the left takes different values when one replaces $\beta$ by $\beta+2\pi i/\log(x)$ while $x^\beta$ does not change. 
    This implies that Riemann's formula should be rewritten using instead the function 
  $$
  \operatorname{Ei}(z)=\int_{-\infty}^z \frac{e^t}{t} d t
  $$
which has a branch cut discontinuity in the complex $z$ plane running from $-\infty$ to $0$. The correct form of \eqref{piriemann} is  
\begin{equation}\label{piriemann1}
f(x)=\operatorname{Li}(x)-\sum_{\Im(\rho)>0}\left(\operatorname{Ei}(\rho \log x)+\operatorname{Ei}(\overline \rho \log x)\right)+\int_x^{\infty} \frac{1}{t^2-1} \frac{d t}{t \log t}-\log 2
 \end{equation}
where  the sum is over the zeros of zeta with positive imaginary part.} the middle term is written as a sum  of the integral logarithm $\operatorname{Li}$ function evaluated on $x^\rho$ where the $\rho$'s are the non-trivial zeros of $\zeta$.\newline 
Now this is obviously a nonsensical sum because $x$ to the power $\rho$ does not change if one replaces $\rho$ by $\rho + 2 \pi i n /\log x$, where $n$ is an arbitrary integer. So the constancy of $\operatorname{Li}(x^\rho)$ on this arithmetic progression implies that there is no possible meaning to the infinite sum because it includes infinitely many repetitions of essentially the same term, since the complex numbers $x^\rho$ accumulate infinitely at least on the circle of radius $x^{1/2}$. In fact,  von Mangoldt gave, many years after your paper, a detailed proof of your theorem, and he was more careful, up to some extent, in writing the formula ... since instead of dealing with $\operatorname{Li}$  he considers a well defined univalent function $\operatorname{Ei}(z)=\operatorname{Li}(\exp(z))$ of the correct variable  $z$ (rather than $\exp(z)$). He then proves your formula replacing the term $\sum \operatorname{Li}(x^\rho)$ by  $\sum \operatorname{Ei}( \rho \log x)$.

 \medskip

  This being said, it is clear from your formula that you knew what is now called the Guinand-Weil explicit formulas. And  the new input  in this formula, but that wouldn't surprise you, is that one considers, instead of the function $1/\log u$ which you use, an arbitrary test function $\phi(u)$ applied on the prime powers, and one writes an equality which involves what is nowadays called the Mellin transform\footnote{The Mellin transform of $\phi(u)$ is defined as $\mathcal{M}(\phi)(s)=\int_0^{\infty} u^{s-1} \phi(u) du$} of this function $\phi(u)$ evaluated on the non-trivial zeros of the zeta function. One obtains in this way an equality which is not really a generalization of your formula in the sense that yours  goes one level of subtlety beyond these classical Guinand-Weil explicit formulas since your test function $1/\log u$ is singular at $u=1$. One thing which is striking in these explicit formulas is that if the test function $\phi(u)$ vanishes outside a finite interval $[1,x]$, then only finitely many primes are involved in its computation, those less  than $x$. But this is of course the case for your formula which involves the function equal to $\phi(u)=1/\log u$ in the interval $[1,x]$ and vanishing outside, so this would not be a surprise to you. 
   
 \medskip
 
But  let me tell you, Master, a fact which I find quite supportive of your hypothesis. The fact is the following. I will  devise a process with which, by using only a few primes, I will be able to reach the first non-trivial zeros of your zeta function with remarkable precision. Now what do I mean by involving a few primes? I will take a concrete situation and I will imagine that I do not know what are the primes beyond 13. So I only know 2, 3, 5, 7, 11 and 13. Now at first it seems that in order to get any grasp on the zeros of zeta, it is necessary to involve all primes, because all of them are involved in the Euler product which is defining this function. And so it seems quite unbelievable that by truncating this Euler product, namely by only considering the terms in the Euler product which are involving the primes smaller than 13, one could get any hint on the zeros of the zeta function.

 \medskip

 So let me describe the process using only notions that are familiar to you. Out of these primes, 2, 3, etc., up to 13, one fabricates a quadratic form. This quadratic form  is similar to the quadratic form you were using when you applied the Dirichlet principle in proving the conformal mapping theorem. It is a quadratic form  on the infinite dimensional space of functions $\phi(u)$ of a positive real variable which vanish outside the interval $[1,13]$.  The value $Q(\phi)$ of the quadratic form  is obtained  by applying the explicit formula to the function $\psi(v)=\int \phi(u)\phi(uv)du/u$. Thus because the function $\phi$ vanish outside the interval $[1,13]$, the function $\psi$ vanishes outside the interval $[1/13,13]$ and one does not need  to use any other prime power than $2,3,4,5,7,8,9,11,13$ to compute $Q(\phi)$. Next, I know how to prove that there is a function $\eta(u)$ realizing the minimum of the quadratic form $Q(\phi)$ while $\int \phi(u)^2du/u=1$. The proof of existence is entirely similar to the proof given by Hilbert in 1900, in his paper "Über das Dirichletsche Prinzip" for the Dirichlet principle. I then take the Mellin transform of the function $\eta(u)$. I also know how to prove\footnote{See \S \ref{sectvW}}
 that  the zeros of this Mellin transform are on the critical line\footnote{Which is normalized here as the imaginary line}. This is proved modulo a  condition of uniqueness of the minimum\footnote{The proof uses a generalization of a Theorem of Caratheodory-Fejer on Toeplitz matrices, obtained in 1911, one needs to assume that the lowest eigenvalue of the quadratic form is simple and even}.

 \medskip

The amazing fact, which I want to report to you, is that when one computes—with the aid of modern computational machines far beyond what was available in your time—the zeros of the Mellin transform of $\eta(u)$,   say its first 50 zeros, one finds an incredible coincidence with the non-trivial zeros of zeta. For instance, for the first zero, you find out that there are 54 decimal places which agree, and so on and so forth. The number of decimal places which agree decreases slowly, and when you reach the 50-th zero, only  few decimal places  agree.

I have computed these differences (upper bound of) between the values computed
using primes less than 13 and the true values of zeros of zeta which I present herewith:
\end{large}
\vspace{0.1cm}

\begin{center}
\rule{0.95\textwidth}{0.4pt}
\end{center}

\begin{center}
\textsc{Differences Between Values (using primes $\le 13$)}
\end{center}

\small
\begin{align*}
&2.60179 \times 10^{-55},\, 4.80071 \times 10^{-52},\, 4.43756 \times 10^{-50},\, 3.89903 \times 10^{-47},\, 7.59453 \times 10^{-46},\, 1.13198 \times 10^{-43},\\
&1.07245 \times 10^{-41},\, 1.2694 \times 10^{-40},\, 4.40141 \times 10^{-38},\, 4.24869 \times 10^{-37},\, 5.86724 \times 10^{-36},\, 3.24443 \times 10^{-34},\\
&2.44517 \times 10^{-32},\, 9.02026 \times 10^{-32},\, 5.13539 \times 10^{-30},\, 7.04142 \times 10^{-29},\, 6.47754 \times 10^{-28},\, 4.96772 \times 10^{-27},\\
& 5.86016 \times 10^{-25},\, 3.76751 \times 10^{-24},\, 1.03779 \times 10^{-23},\, 3.52722 \times 10^{-22},\, 3.03977 \times 10^{-21},\, 5.66201 \times 10^{-20},\\
&1.41755 \times 10^{-19},\, 2.19821 \times 10^{-18},\, 6.31599 \times 10^{-17},\, 1.42037 \times 10^{-16},\, 4.34328 \times 10^{-16},\, 4.47113 \times 10^{-15},\\
& 7.01522 \times 10^{-14},\, 3.81989 \times 10^{-13}, \, 5.99581 \times 10^{-13},\, 4.26414 \times 10^{-11},\, 1.10653 \times 10^{-10},\, 1.95651 \times 10^{-10},\\
&5.20728 \times 10^{-10},\, 2.05031 \times 10^{-9},\, 3.42274 \times 10^{-8},\, 2.10931 \times 10^{-7},\, 2.23714 \times 10^{-7},\, 5.95608 \times 10^{-7},\\
& 5.77737 \times 10^{-6},\, 0.000141389,\, 0.000556111,\, 0.000720794,\, 0.000314865,\, 0.0209081,\, 0.00313565,\, 0.00212727
\end{align*}
\normalsize

\begin{center}
\rule{.95\textwidth}{0.4pt}
\end{center}\vspace{0.1cm}
\begin{large}
 What this is saying is that we have a firm grasp on your zeros, without at any point involving the infinity of the collection of all primes. And moreover we know a priori that  all zeros of the Mellin transform of $\eta(u)$ are on the critical line.

 \medskip

 What we do not know is that, when we increase the upper limit, which was  $x=13$ here, the corresponding set of zeros will converge towards the  zeros of zeta. This is something which at this point is not proved. On the other hand, it seems that the abstract reason why your conjecture is true is now within reach since we know that the zeros of the Mellin transforms of the minimal eigenvectors $\eta_x$ are purely imaginary, and we expect that the $\eta_x$ converge to the function whose Mellin transform is your function $\Xi(it)$. The result would then follow from Hurwitz's theorem, which implies that all the zeros of the limit of a convergent sequence\footnote{uniformly on compact subsets} of holomorphic functions whose zeros are on a fixed line are still on that line.
 
 \vspace{.1in}

\centerline{With respect and admiration, }

\vspace{.1in}

\centerline{Alain Connes}

\end{large}

\endgroup
 
 \newpage
 
\section{The strategy and the  next small steps}
We first explain in \S \ref{sectvW} the general result which provides a wealth of entire functions all of whose zeros are on the real line.
The goal then is to try to show that, when the upper bound $x$ (which was $x=13$ above) tends to infinity, the minimal eigenvectors $\eta_x$ converge in a suitable sense so that their Fourier transforms, for the duality between $\R^*_+$ and $\R$, will converge uniformly on compact subsets to the holomorphic function that Riemann had introduced as the $\Xi$-function in his paper. In fact we shall follow  his footsteps and consider the explicit formula he gave for the function $k$ whose Fourier  transform is  $\Xi(t)$. We shall  rewrite in \S \ref{sectxi} this function $k$ as the image by a summation map $\cE$ of a linear combination $h$ of two Hermite functions $h_0,h_4$. The issue, then, is to show the convergence of the minimal eigenvectors $\eta_x$ to the function $\cE(h)$. In fact, one first needs to recenter $\eta_x$ to $\theta_x(u)=\eta_x(x^{1/2}u)$ which has its support  in $[x^{-1/2},x^{1/2}]$. The plan, in order to achieve this goal, is to break it into several small steps. In the second step, \S \ref{sectprolate}, one introduces, given $\lambda>1$, the prolate spheroidal wave functions which give an adaptation $h_{n,\lambda}$ of the hermite functions $h_n$ to functions with support in the interval $[-\lambda,\lambda]\subset \R$.  The functions $h_{n,\lambda}$ are eigenfunctions of a modification of the Hermite operator \eqref{herm}, the prolate wave operator $\pw_\lambda$ of \eqref{wop}, obtained by adding to the Hermite operator $\mathbf{H}$ the only multiple of the  square of the scaling so that the sum commutes with the orthogonal projection on functions with support in the interval $[-\lambda,\lambda]\subset \R$. By construction $\pw_\lambda$ admits two regular singular points at the boundary $\pm \lambda$ of the interval and an irregular singularity at $\infty$. The third step \S \ref{poisson}, consists in making an educated guess for an approximation  of the minimal eigenvector $\theta_x$ by using the same process as in the simple case of \S \ref{sectxi}. Thus one simply replaces the Hermite functions $h_0,h_4$ by their localizations $h_{0,\lambda},h_{4,\lambda}$, for $\lambda=x^{1/2}$,
and one applies the summation map $\cE$ to the  linear combination of these two functions with vanishing integral. One obtains in this manner a function $k_\lambda$ which approximates $\theta_x$ on the interval $[x^{-1/2},x^{1/2}]$. In \S \ref{sectconv} one shows that, when $\lambda\to \infty$, the Fourier transforms $\widehat k_\lambda$ converge towards the $\Xi$ function of Riemann, uniformly on closed substrips of the open  strip of width $1$ neighborhood of the real line. Finally in \S \ref{miss} we point to the missing final steps.
 \subsection{The Fourier transform of $\theta_x$ has all its zeros on the real line}\label{sectvW}
The proof of this result follows from a theorem shown in a joint paper \cite{CS} with Walter van Suijlekom, entitled "Quadratic Forms, Real Zeros  and Echoes of the Spectral Action". The precise statement is :
\begin{thm}\label{mainintro}
Let \( L > 0 \), \( \mathcal{D} \) be a real  distribution on the interval \( [0, L] \) and $\tilde \cD$ the associated even distribution on \( [-L, L] \). Assume that the quadratic form with Schwartz kernel \( \tilde \cD(x - y) \) defines a lower-bounded selfadjoint operator  on \( L^2([-\frac{L}{2}, \frac{L}{2}]) \), and that the minimum of its spectrum is a simple, isolated eigenvalue, with even eigenfunction $\eta$. Then all the zeros of the entire function $\widehat \eta(z)$, $z\in \C$, Fourier transform  of $\eta$ lie on the real line.
\end{thm}
The proof of this theorem is based on the special form of the matrix of the quadratic form
 in the trigonometric orthonormal basis, the construction for finite matrices of that special form of a selfadjoint operator and the above mentioned Hurwitz theorem to pass to the limit when the size of the matrices tends to infinity. The validity of the result on the zeros of the Fourier transform for finite truncations \cite{CS} plays a key role in numerical computations\footnote{It is used for the trigonometric truncation at $N=100$ in the computation presented in the letter} and allows one to approximate the zeros of $\widehat \eta(z)$ by the spectrum of a rank one perturbation of the periodic Dirac operator, obtained using the Dirichlet kernel analysed by Riemann in \cite{Riemann1}.

\subsection{Riemann's $\Xi$ function and Hermite functions}\label{sectxi}
Riemann adopts the following notation
$$
\psi(x):=\sum_1^\infty e^{-n^2\pi x}
$$
and then writes, with $\xi(s):=\frac{1}{2} s(s-1) \pi^{-s / 2} \zeta(s) \Gamma\left(\frac{s}{2}\right)$, $\Xi(t):=\xi(s)$ for $s=\frac 12 +it$, the equality\footnote{The traditional notation for this function is $\Xi (t)$ instead of Riemann's $\xi(t)$.}
$$
\Xi(t)=4 \int_1^{\infty} \frac{d\left(x^{\frac{3}{2}} \psi^{\prime}(x)\right)}{d x} x^{-\frac{1}{4}} \cos \left(\frac{1}{2} t \log x\right) d x
$$
He then obtains, using the equality $k(u)=k(u^{-1})$ from the Poisson formula, where
$$
 k(u):= u^{1/2}\frac \pi 2\sum_1^\infty n^2u^2\left(2 \pi  n^2 u^2-3\right)e^{-\pi n^2u^2}
$$ 
that  $\Xi(t)$ is the Fourier transform  of the function  $k(u)$ :
$$
\Xi(t)=2\int_1^\infty k(u)\cos \left( t \log u\right)d^*u=\int_0^\infty k(u)u^{it}d^*u, \ \
$$ 
 Moreover introducing the notation 
\begin{equation}\label{emap}	
\cE(f)(u):=u^{1/2}\sum_1^\infty f(nu)
\end{equation}
 one has 
\begin{equation}\label{ku}	
k(u)=\cE(h)(u), \ \ h(u)=\frac \pi 2 u^2\left(2 \pi   u^2-3\right)e^{-\pi u^2}.
\end{equation}
The function $h(u)$ can be characterized as follows. One considers the Hermite operator (harmonic oscillator)
 \begin{equation}\label{herm}	\mathbf{H}f(u):=-f''(u)+4 \pi ^2 u^2 f(u)
 \end{equation} 
 and lets  $h_n$ be the normalized eigenfunction for the eigenvalue $2\pi(1+2n)$. These functions are even for $n$ even and invariant under the Fourier transform for $n$ multiple of $4$.
\begin{fact}\label{hermfact}
The $\Xi$ function of Riemann is the Fourier transform of $k=\cE(h)$ where 	$h$ is, up to a multiplicative scalar, the only linear combination of $h_0, h_4$ with vanishing integral\endnote{More precisely one has 
\begin{equation}\label{xih}	
  h=\frac{\sqrt{3}}{4\ 2^{3/4}}h_4-\frac{3}{2^{17/4}}h_0,  \ \  \Vert h_j\Vert=1, \  \  \Vert h\Vert=\frac{\sqrt{33}}{2^{17/4}}.
\end{equation}}.
\end{fact}

\subsection{The prolate wave functions enter the scene}\label{sectprolate}

 The next step requires introducing a family of functions whose pivotal role in signal transmission was uncovered by Slepian, Pollak, and Landau at Bell Labs during the 1960s. Their work addressed a fundamental question first posed by Claude Shannon in his foundational study on entropy\footnote{C.~E.~Shannon, \emph{Communication in the Presence of Noise}, Proc.\ IRE \textbf{37} (1949), 10--21.}: \emph{To what extent can functions that are band-limited also be time-limited?} Shannon had estimated that the number 
 \( N \) of orthogonal signals that can be packed into a time window \( T \) and frequency band \( W \) satisfies \( N \simeq 2TW \).
 Building on this, Slepian and his collaborators sought to maximize the information transmittable over a limited time within a given bandwidth — a problem central to efficient communication systems.

To describe their result, fix the time and frequency interval as \( [-\lambda, \lambda] \subset \mathbb{R} \), and let \( P_\lambda \) be the projection in the Hilbert space \( L^2(\mathbb{R}) \) (square-integrable  functions) defined by multiplication with the characteristic function of the interval \( [-\lambda,\lambda] \). Let \( \widehat{P_\lambda} := \mathbb{F}_{e_{\mathbb{R}}} P_\lambda \mathbb{F}_{e_{\mathbb{R}}}^{-1} \) be its conjugate under the Fourier transform,\footnote{The normalized Fourier transform is defined by
\[
\mathbb{F}_{e_{\mathbb{R}}}(f)(y) := \int_{\mathbb{R}} f(x) e^{2\pi i x y} \, dx.
\]
Note that on even functions one has \( \mathbb{F}_{e_{\mathbb{R}}}^{-1} = \mathbb{F}_{e_{\mathbb{R}}} \).}

In their seminal paper, \cite{Slepian0},
%\footnote{D.~Slepian and H.~Pollak, \emph{Prolate Spheroidal Wave Functions, Fourier Analysis and Uncertainty}, The Bell System Technical Journal (1961), 43--63.} 
Slepian and Pollak showed that the compressed Fourier transform \( P_\lambda \mathbb{F}_{e_{\mathbb{R}}} P_\lambda \) commutes with the second-order differential operator
\begin{equation}\label{wop}
\pw_\lambda := -\partial_x\left( (\lambda^2 - x^2) \partial_x \right) + (2\pi \lambda x)^2.
\end{equation}

It follows that this operator commutes on $L^2(\R)^{\mathrm{ev}}$ (even functions)  with
\[
P_\lambda \widehat{P_\lambda} P_\lambda = P_\lambda \mathbb{F}_{e_{\mathbb{R}}} P_\lambda \mathbb{F}_{e_{\mathbb{R}}} P_\lambda = \left(P_\lambda \mathbb{F}_{e_{\mathbb{R}}} P_\lambda\right)^2,
\]
which enables one to diagonalize the \emph{angle operator}\endnote{The notion of angle between two closed subspaces of a Hilbert space — or equivalently, between two orthogonal projections \( P \) and \( Q \) — was studied by J.~Dixmier in: \emph{Position relative de deux variétés linéaires fermées dans un espace de Hilbert}, Revue Sci.\ \textbf{86} (1948), 387--399. The key idea is that a pair \( (P, Q) \) corresponds to a unitary representation of the infinite dihedral group. Its irreducible representations are classified by an angle \( \alpha \), which in the real plane represents the angle between two real lines. The relation \( P\cos^2 \alpha = P Q P \) determines the angle operator on each irreducible component.} \( \alpha_\lambda \) between the projections \( P_\lambda \) and \( \widehat{P_\lambda} \) defined by the identity
\[
P_\lambda \cos^2(\alpha_\lambda) = P_\lambda \widehat{P_\lambda} P_\lambda.
\]
\begin{fact} The eigenvalues of the operator $P_\lambda \widehat{P_\lambda} P_\lambda$ in $L^2([-\lambda, \lambda])^{\mathrm{ev}}$ are simple and form a decreasing sequence $\nu_n(\lambda)$, $n\geq 0$, $\nu_n(\lambda)\to 0$ for $n\to \infty$, such that $1>\nu_0(\lambda)>\nu_1(\lambda>\ldots>0$. The corresponding eigenfunctions are the \emph{prolate spheroidal wave functions} of even index $h_{2n,\lambda}$ where $h_{m,\lambda}$ is the $m+1$-th  eigenfunction of the prolate wave operator \eqref{wop} in $L^2([-\lambda, \lambda])$.	
\end{fact}

The operator \( \alpha_\lambda \) has a finite number\footnote{With our normalization, time duration is \( T = 2\lambda \) and bandlimit is \( W = \lambda \).} \( \sim 4\lambda^2 \) of small nonzero eigenvalues.  Following the standard notation in the literature, the prolate spheroidal wave functions are denoted
$$
h_{n, \lambda}(x):=\mathrm{PS}_{n, 0}\left(2 \pi \lambda^2, \frac{x}{\lambda}\right)
$$
where the subscript $n$ indexes the eigenfunction, the second subscript $0$ indicates the angular order-a remnant of the operator's origin in separating variables for the Laplacian on prolate spheroids-and the first argument $2 \pi \lambda^2$ is the dimensionless bandwidth parameter.
Each \( h_{n,\lambda} \) is supported in \( [-\lambda,\lambda] \) and extended by zero outside. The function \( h_{n,\lambda} \) is even when \( n \) is even. The Fourier transform of \( h_{2m,\lambda} \), restricted to \( [-\lambda, \lambda] \), satisfies
\[
\mathbb{F}_{e_{\mathbb{R}}}(h_{2m,\lambda}) = \chi_m h_{2m,\lambda},
\]
where $\chi_m^2=\nu_{m}$ and the sign of $\chi_m$ is $(-1)^m$. For our purpose we need the eigenvalues $\chi_0, \chi_2$ corresponding to the eigenfunctions $h_{0,\lambda},h_{4,\lambda}$ and both $(1-\chi_0)$, and $(1-\chi_2)$ tend exponentially to $0$ as functions of $x=\lambda^2$.
%\begin{figure}[H]
%\begin{center}
%\includegraphics[scale=0.35]{legend1.png}\\
%\caption{Graphs of $\log(1-\chi_0)$ and $\log(1-\chi_2)$ as functions of $x=\lambda^2$.\label{fpro}}
%\end{center}
%\end{figure}

\subsection{The Poisson formula and the approximation $k_\lambda$ of $\theta_x$, $\lambda=\sqrt x$.}\label{poisson}
\begin{figure}[H]
\begin{center}
\includegraphics[scale=0.35]{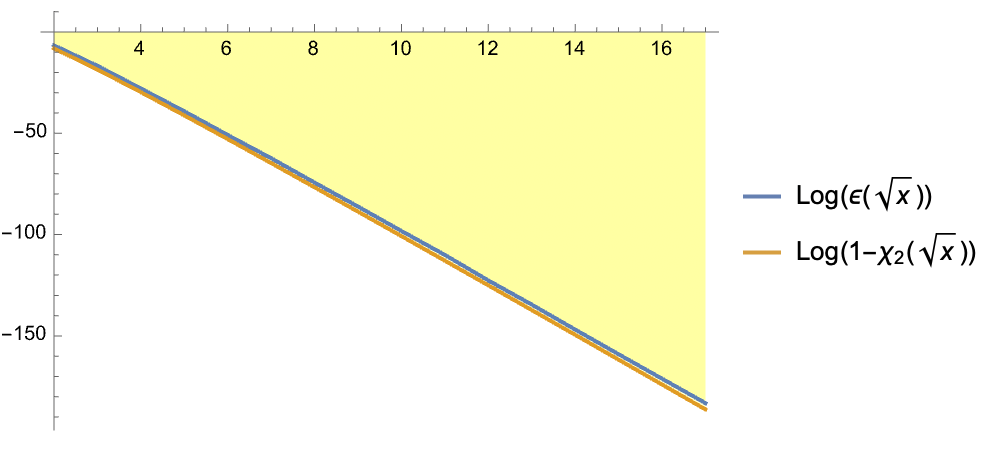}\\
\caption{Graphs of $\log(\epsilon(\sqrt x)))$ and $\log(1-\chi_2(\sqrt x)))$ as functions of $x$.\label{fpro1}}
\end{center}
\end{figure}

Let $\lambda>1$, and $QW_\lambda$ be the restriction of the Weil quadratic form to test functions whose support is within the interval  $[\lambda^{-1},\lambda]$.  By the result of André Weil discussed in \S \ref{sectweilpos}, the positivity of $QW_\lambda$ for all $\lambda>1$ is equivalent to RH. This positivity can be proved for small values of $\lambda$ (see \cite{yoshida}, \cite{weilpos}). There is (see \cite{VJ,c2m2b}) for each $\lambda>1$ a canonical lower bounded, unbounded selfadjoint  operator  $A_\lambda$ with compact resolvent, in the Hilbert space  $L^2\left(\left[\lambda^{-1}, \lambda\right], d u/u\right)$ such that 
\begin{equation}\label{qtoa}
	Q W_\lambda(f, f)=\langle A_\lambda f\mid f\rangle
\end{equation}
  The numerical computation of the smallest eigenvalue $\epsilon(\lambda)$ of $A_\lambda$, done in \cite{VJ}, shows that $\epsilon(\lambda)$ tends exponentially fast to $0$ as a function of $\mu=\lambda^2$. In fact a careful analysis reveals  a striking similarity (Figure \ref{fpro1}) between the behavior of $\epsilon(\lambda)$ and of the angular function $1-\chi_2(\lambda)$.\newline In terms of the length $L=2\log \lambda$  of the support $[\lambda^{-1}, \lambda]$ of the test functions for $Q W_\lambda$, the convergence to $0$ of the minuscule  quantities like $1-\chi_2$ is {\it exponential of exponential}\footnote{see \cite{Fuchs}, Theorem 1. Note  that  $\chi_k(\lambda)^2=\lambda_{2k}(a)$ with $a=\sqrt{2\pi}\,\lambda$ in the notations of this theorem }: $$1-\chi_2\sim\frac{2^{14}}{3} \sqrt{2} \pi ^5 e^{-4 \pi  e^L+9/2 L} $$.  \newline 
    In \cite{VJ}, we gave a construction  numerically\endnote{This fact is verified in \cite{VJ},  where we compare the eigenvectors of the Weil quadratic form $QW_\lambda$  associated to its smallest eigenvalues with the orthogonalisation of the prolate vectors  obtained using the technique outlined above, from the prolate spheroidal wave functions.} justified of the  eigenvectors associated to the first minuscule eigenvalues of $Q W_\lambda$, using prolate spheroidal wave functions associated to the interval $[-\lambda,\lambda]$. In particular this gives an educated guess for an approximation of the eigenvector associated to the smallest eigenvalue $\epsilon(\lambda)$ of $A_\lambda$. In  agreement with Fact \ref{hermfact}, it is \begin{equation}\label{ktoh}
  k_\lambda(u):=\cE(h_\lambda)(u), \ \ \forall u\in [\lambda^{-1},\lambda]
  \end{equation}
   where $h_\lambda$ is, up to a multiplicative scalar, the only linear combination of $h_{0,\lambda}, h_{4,\lambda}$ with vanishing integral\footnote{Note that the computation of $\cE(h_\lambda)(u)$ for $u\in [\lambda^{-1},\lambda]$ only involes the sum over the integers $\leq \lambda^2$}. 
   
   The conceptual justification of this formula is as follows: The range of the map $\cE$ is contained in the radical of the global Weil quadratic form (see \cite{CMbook}), but   RH implies that $QW_\lambda$ is strictly positive and that  its radical is $\{0\}$ so   we should expect that the domain of $QW_\lambda$ cannot contain any non-zero element of the range of the map $\cE$. One can nevertheless  construct   functions  with  support in $[\lambda^{-1},\lambda]$ which are in the ``near radical'' of the Weil quadratic form as follows.   If the support of the even function $f$ is contained in the interval $[-\lambda,\lambda]\subset \R$, the support of $\cE(f)$ is contained in $(0,\lambda]\subset \R_+^*$. On the other hand, the Poisson formula, using the conditions $f(0)=\widehat f(0)=0$ to define the codimension $2$ subspace $\mathcal S_0\subset \mathcal S(\R)$ of the Schwartz space, gives 
\begin{equation}
\cE(\widehat f)(x)=\cE(f)(x^{-1}) \qqq f\in \sr0\label{poissonintro}
\end{equation}
 which shows that the support of $\cE(f)$ is contained in $[\lambda^{-1},\infty)$ provided the support of the even function $\widehat f$ is contained in the interval $[-\lambda,\lambda]\subset \R$. The obstruction to obtain an element $\cE(f)$ of the radical of $QW_\lambda$ is the equality $P_\lambda \cap \widehat P_\lambda=\{0\} $, where  $P_\lambda$ and $\widehat{P_\lambda}$ are as above. But as explained in \S \ref{sectprolate} these two projections nearly intersect, and after taking care of the two conditions $f(0)=\widehat f(0)=0$,  the restriction of $\cE(f)$ to the interval $[\lambda^{-1},\lambda]$ gives rise to the function $k_\lambda$ of \eqref{ktoh}  on which $QW_\lambda$ takes non-zero, but extremely small values, thus giving a natural guess as  an approximation of the eigenvector associated to the smallest eigenvalue $\epsilon(\lambda)$ of $A_\lambda$.

 \subsection{Convergence of the Fourier tranforms  $\widehat k_\lambda\to \widehat k$}\label{sectconv}
 Thanks to the classical estimates on the convergence of the prolate wave functions towards the Hermite-Weber function one controls the convergence of $k_\lambda$ of \eqref{ktoh} towards $k=\cE(h)$ and this gives
 \begin{fact}\label{hermfact1} The Fourier transform of $k_\lambda$ converges, when $\lambda\to \infty$, towards the $\Xi$-function of Riemann uniformly on closed substrips of the open  strip $\Im(z)<\frac 12$.
\end{fact}
The difference is controlled on the line $\Im(z)=\alpha$ where $\alpha\in (-\frac 12, \frac 12)$ by $c\lambda ^{-\frac{1}{2}-\alpha }(1-2 \alpha)^{-1}$ where $c$ is a finite constant.

\subsection{Remaining steps}\label{miss} In order to apply Theorem \ref{mainintro} one needs to show that the smallest eigenvalue of the Weil quadratic form $QW_\lambda$ is simple with even eigenvector. The analogue of this property is known for the prolate wave operator. Moreover it still remains to show that  $k_\lambda$ is a sufficiently good approximation of $\theta_x$, $\lambda=\sqrt x$.

\section{Geometric Perspectives}
The geometric approach pursued in the joint work  \cite{ConnesConsaniTopos}  described in \cite{Essay}  unveiled a novel geometric landscape which is a natural analogue, for the global field $\Q$, of the geometry  associated to  global fields in finite characteristic. In particular the Frobenius correspondences make sense on the square of the scaling site and the completed Riemann zeta function appears as the analogue of the Hasse-Weil generating function. The  handling of the delicate principal values involved in the Riemann-Weil explicit formula required to use the adelic incarnation of the scaling site which we now understand as its class field theory counterpart \cite{ccknots}.\newline
In the paper \cite{weilpos}, we began to exploit the adelic interpretation of the Frobenius correspondences for the number field $\Q$ in order to gradually investigate Weil's positivity. As explained in \S \ref{sectweilpos} the Weil positivity, which only involves finitely many primes at a time, is equivalent to RH. In the adelic geometric framework, the correspondences are encoded by Schwartz kernels, which are distributions in the square of the adelic data that play the role of the class field theory counterpart of the geometric curve. More generally the connection between the operator theoretic and the geometric viewpoints is effected by the Schwartz kernels associated to operators. By implementing the additive structure of the adeles, one sees that the Schwartz kernel of the scaling operator corresponds geometrically to the divisor of the Frobenius correspondence.
 \begin{graybox} The additive structure of the adeles of $\mathbb{Q}$ allows one to write the Schwartz kernel $k(x, y)$ of the scaling action $f(x) \mapsto f(\lambda x), \lambda \in \mathbb{R}_{+}^*$, in the form (with $\underline{\delta}$  the Dirac distribution)
$$
k(x, y)=\underline{\delta}(\lambda x-y). 
$$
\end{graybox}
\subsection{Archimedean trace formula} 
There is a very important  parallel between the Weil quadratic form and the trace formula, on one hand, and the  world of information theory. The starting point of this parallel is to rewrite the archimedean case of the  trace formula of \cite{Co-zeta} by involving now two independent parameters which play the role of the time limitation and frequency limitation in the approach of Shannon, Slepian and their  collaborators concerning the transmission of information. One works in the Hilbert space $L^2(\mathbb{R})_{\mathrm{ev}}$ of even square integrable functions. The scaling action of $\mathbb{R}_{+}^*$ is defined as  
$$(\vartheta_\lambda \xi)(v):=\lambda^{-1 / 2} \xi\left(\lambda^{-1} v\right), \ \vartheta(f)=\int f(\lambda) \ \vartheta_\lambda d^*\lambda, \quad \forall \xi \in L^2(\mathbb{R})_{\mathrm{ev}}$$
The parameter $T$ which gives the limitation in time to the interval $[-T,T]$ defines the orthogonal projection $P_T$ and the parameter $W$ which gives the limitation in frequency to the interval $[-W,W]$ defines the projection $\widehat P_W$. One can write the archimedean case of the  trace formula of \cite{Co-zeta} as follows, where $W_\infty:=-W_\R$ was defined in \eqref{arch}:
 \begin{graybox}
 \begin{equation}\label{weilq}
 		W_\infty(f)= \log( TW)\ f(1)+\operatorname{Trace}\left(\vartheta(f)\left(1-P_T-\widehat P_W\right)		\right)
 			\end{equation}
 			\end{graybox}   
We view this formula as a bridge between the explicit formulas and the world of information theory where the projections $P_T$ and $\widehat P_W$ play a central role. It will be extended below in \S \ref{sectsemiloc} to incorporate the contribution of the primes to the explicit formulas. 

\subsection{Archimedean  Weil positivity}
The key ingredient is the semilocal trace formula, which in our paper was used in the simple case when no primes are involved. What we found is that not only Weil's positivity holds in that case, as explained in \S \ref{sectweilpos}, but the main source of positivity is due to the Sonin space, first introduced in the context of RH by Burnol \cite{burnol-1,burnol-2,burnol-3}. The Sonin space  ${\mathfrak S}_\lambda$ is the space of square integrable functions which vanish identically in the interval $[-\lambda,\lambda]$ as well as their Fourier transform. It is by construction the orthogonal of the ranges of the projections $P_\lambda$ and $\widehat{P_\lambda} $. We denote by ${\mathfrak S}$ the orthogonal projection in 
$L^2(\mathbb{R})_{\mathrm{ev}}$ to the Sonin space ${\mathfrak S}_1$ (for $\lambda=1$). One then has \cite{weilpos},
\begin{thm} Let $g \in C_c^{\infty}\left(\mathbb{R}_{+}^*\right)$ have support in the interval $\left[2^{-1 / 2}, 2^{1 / 2}\right]$ and Fourier transform vanishing at $\frac{i}{2}$ and 0 . Then the following inequality holds
$$
W_{\infty}\left(g * g^*\right) \geq \operatorname{Tr}\left(\vartheta(g) {\mathfrak S} \vartheta(g)^*\right)
$$
\end{thm}
%\subsection{The semilocal set-up and trace formula}
\subsection{The semilocal adele class space}
Underlying the semilocal trace formula discussed in \S \ref{sectsemiloc}, are the semilocal adele class spaces.
These geometric spaces $Y_S$ are associated with a finite set $S$ of places of $\Q$ containing the infinite place. By construction $Y_S$ is the quotient
\begin{equation}\label{YQS}
Y_{S}:=\A_{S}/\Gamma_S,\ \  \A_{S}=\prod_{v\in S} \Q_v
\end{equation} 
of the adelic space 
 		 		 product of the local fields completions of the global field $\Q$ at the places $v\in S$. The group $\Gamma_S$ is the subgroup of $\Q^\times$ defined by  
 \begin{equation}\label{GL1QS}
\Gamma_S:= \{ \pm p_1^{n_1} \cdots p_k^{n_k} \, :\,  p_j
\in S \setminus\{ \infty \} \,,\, n_j\in \Z\}\subset \Q^\times
\end{equation}
The ring $ \A_{S}$ contains $\Q$ as a subring using the diagonal embedding and this induces the action of $\Gamma_S$ on $ \A_{S}$ by multiplication. The  semilocal adele class spaces are best encoded by the noncommutative algebras $\cS(\A_S)\rtimes \Gamma_S$ cross products of the Bruhat-Schwartz algebras $\cS(\A_S)$ of functions on semilocal adeles, by the multiplicative groups $\Gamma_S$.  These noncommutative algebras constitute a sheaf of algebras over $\Spec\,\Z$. The group $ \Gamma_S$ is properly understood as sections $\Z_S^\times$ of the sheaf $\mathbf{G}_m$ on the open set complement of $S$ in $\Spec\,\Z$. One shows that the semilocal Bruhat-Schwartz algebra $\cS(\A_S)$ form a sheaf $\sheaf$ of commutative algebras over $\Spec\,\Z$. One then obtains the following result which establishes the compatibility of the noncommutative geometric constructions with the algebraic geometry of $\Spec\,\Z$. 
\begin{thm}	
\label{structure1intro}  
 		\begin{enumerate}
 			\item The algebraic cross product $\sheaf\rtimes \mathbf{G}_m$ defines a sheaf  of algebras on  $\Spec\Z$ such that for every finite set of places $S\ni \infty$
 			$$
 			\left(\sheaf\rtimes \mathbf{G}_m\right)(S^c)=\cS(\A_S)\rtimes \Z_S^\times 			$$
 			\item The stalk of $\sheaf\rtimes \mathbf{G}_m$ at the generic point  is the global cross product $\cS(\A_\Q)\rtimes \Q^\times$. 			\item The global sections of $\sheaf\rtimes \mathbf{G}_m$ form the cross product $\cS(\R)\rtimes \{\pm 1\}$.
 		\end{enumerate}
\end{thm}
The function spaces involved in the semilocal trace formula of \S \ref{sectsemiloc} are best understood conceptually as the Hochschild homology of the semilocal algebras.

\subsection{The semilocal trace formula}\label{sectsemiloc}
The remaining  difficulty in proving that the eigenvectors $\theta_x$ converge to the function $k=\cE(h)$ of Fact  \ref{hermfact} is to effectively compare $\theta_x$ with $k_\lambda$ for $\lambda=\sqrt x$. The numerical evidence was shown in \cite{VJ} where the comparison was extended to the eigenvectors of $QW_\lambda$ corresponding to the first minuscule eigenvalues, using the Gram-Schmidt orthogonalisation of vectors of the form $\cE(\psi)$ where the $\psi$ are constructed using the next prolate wave functions.\newline
  As a step towards a conceptual justification of this numerical fact one has the semilocal trace formula of \cite{Co-zeta}. It provides a trace formula   representation of the Weil quadratic form $QW_\lambda$ which is perfectly analogous to \eqref{weilq}, but now gives the contribution of the primes $p\in S$ to the explicit formula. It takes the form 
   \begin{graybox}
   \begin{equation}\label{weilqb}
 		-\sum_{v \in S} W_v(f)=\log( TW)\ f(1)+\operatorname{Trace}\left(\vartheta(f)\left(1-P_T^S-\widehat P_W^S\right)		\right)
 			\end{equation} 
 			 \end{graybox} 
 			where the projections $P_T^S$ and $\widehat P_W^S$ are defined as in the archimedean case using the module\endnote{ 
The module extends to a multiplicative map $\vert \bullet \vert_S$ from the ring $\A_S=\prod_{v\in S} \Q_v$ to $\R_+$, and by construction this map passes to the quotient as a map ${\rm Mod}_S:Y_S=\A_{S}/\Gamma\to \R_+$.
\begin{equation}\label{module}
{\rm Mod}_S(u):=\vert (u_v)_{v\in S} \vert_S=\prod \vert u\vert_v   \in \R_+
\end{equation}. 
The idele and idele class groups
\begin{equation}\label{GL1AS}
\GL_1(\A_{S})= \prod_{p\in S} \GL_1(\Q_p), \ \ C_{S}=\GL_1(\A_{S})/\Gamma
\end{equation}
act naturally by multiplication on the quotient $Y_{S}$ and the orbit of
$1\in \A_{S}$  gives an embedding
$
C_{S}\to Y_{S}.
$
The complement of $C_S$ in $Y_S$ is of measure zero for the product of the  Haar measures of the additive groups of the local fields (which is preserved by the action of the countable group $\Gamma$). Using the Radon-Nikodym derivative of the Haar measures of the multiplicative groups with respect to the Haar measure of the additive groups, one obtains a unitary identification
\begin{equation}\label{wS}
 w_S :L^2(Y_{S})\to L^2(C_{S})
 \end{equation} (see \cite{CMbook} Proposition 2.30).
 We also recall (see  \cite[Eqs. (2.223) and (2.239)]{CMbook} that  $C_{S}$  is a
 modulated locally compact group with module 
\begin{equation}\label{Modulus}
{\rm Mod}_S(\lambda) =   |\lambda|_S :=  \prod_{p\in S} |\lambda_p |, \quad
 \forall \lambda = (\lambda_p) \in  C_{S} 
 \end{equation}
which is (non-canonically) isomorphic to \,$\R^*_+ \times K_S $ ,
where $K_S$ is the kernel of  ${\rm Mod}_S$. The trace formula in the form of \eqref{weilqb} is the specialization of the general trace formula of \cite{Co-zeta} to the $K_S$ invariant part.}.
 	\subsection{The infrared and ultraviolet regimes} \label{ultraviolet}
 	In the framework of noncommutative geometry, the encoding of a geometric space by a spectral triple $(\mathcal{A}, \mathcal{H}, D)$ reveals information through two complementary spectral regimes of the Dirac operator D:

	•	The ultraviolet (UV) regime corresponds to the high-energy part of the spectrum, i.e., the behavior of the eigenvalues of $D$ at infinity. This regime captures the infinitesimal structure of the space and determines its local geometric invariants, such as dimension, volume form, and curvature. In physics, this is akin to probing a system at very short distances or high momenta. Mathematically, the UV behavior governs the asymptotics of the heat kernel and enters crucially in the spectral action principle.
	
	•	The infrared (IR) regime concerns the low-lying part of the spectrum of D, particularly the small eigenvalues. This part reflects the global topological and geometric features of the space, such as connectedness, volume growth, and index-theoretic quantities. In physics, IR behavior refers to low-energy phenomena and large-scale properties. It plays a central role in understanding long-range correlations and topological effects.

 In particular, such a dual perspective allows noncommutative geometry to access both local spectral invariants and global arithmetic or topological features through a unified operator-theoretic language.\newline 
 For the infrared regime we construct in \cite{c2mzeta} self‑adjoint operators \( \dln \) obtained as rank‑one perturbations of the spectral triple associated with the scaling operator on the interval \([ \lambda^{-1}, \lambda ]\) and whose spectrum coincides with the stunning approximation of the low lying zeta zeros as described in the letter to Riemann. We further compute the regularized determinants \( \det\nolimits_{\mathrm{reg}}(\dln - z) \) of these operators and discuss the analytic role they play in controlling and potentially proving the above result by showing that, suitably normalized, they converge towards the Riemann $\Xi$ function.\newline  
 For the ultraviolet regime, we shall briefly describe  in Section \ref{cmpro} the results of \cite{CM} which show that the prolate wave operator provides a self-adjoint operator that matches the ultraviolet behavior of the zeros. As a preparation one can use  the explicit formulas  to compute the heat expansion, assuming RH, of an operator whose spectrum is formed of the imaginary parts of non-trivial zeros of $\zeta(z)$ :
 \begin{thm}\label{heatkernelintro}\cite{heatc} Assume RH and let $D$ be the self-adjoint operator whose spectrum is formed of the imaginary parts of non-trivial zeros of the Riemann zeta function. One then has the asymptotic expansion for $t\to 0$
	\begin{equation}\label{asymp}
	\Tr(\exp(-tD^2))\sim \frac{\log \left(\frac{1}{ t}\right) }{4 \sqrt{\pi } \sqrt{t}}-\frac{(\log 4\pi +\frac 12\gamma)}{2\sqrt{ \pi } \sqrt{t}}+2\exp(t/4)+\sum a_n t^{n/2}
\end{equation}
where $a_0=-\frac 14$ and for $k>0$,  using Bernouilli numbers $B_j$ and Euler numbers $E(k)$,
$$a_{2k-1}=
  \frac{\Gamma(k) \left(2^{2k-1}-1\right) B_{2 k}}{2\sqrt \pi(2 k)!}, \ \ a_{2k}=-\frac 14 \, \Gamma(k+\frac 12)\,\frac{E(2k)}{\sqrt \pi(2k)!}.
  $$
\end{thm}
The Euler numbers are defined as 
\begin{equation}\label{euler}
E(2n):=\sum _{k=1}^{2 n} \left(-\frac{1}{2}\right)^k \sum _{j=0}^{2 k} (-1)^j \binom{2 k}{j} (k-j)^{2 n}
 \end{equation}
 One has the asymptotic behavior when $k\to \infty$
$$
\frac{E(2k)}{(2k)!}\sim (-1)^k 2^{2k}\frac 4 \pi \pi^{-2k}
$$
which shows that the asymptotic expansion \eqref{asymp} is by no means convergent since its general coefficient $a_n$ diverges like a factorial.\newline
 
\subsection{The prolate wave operator}\label{cmpro}
\begin{figure}[H]	\begin{center}
\includegraphics[scale=1]{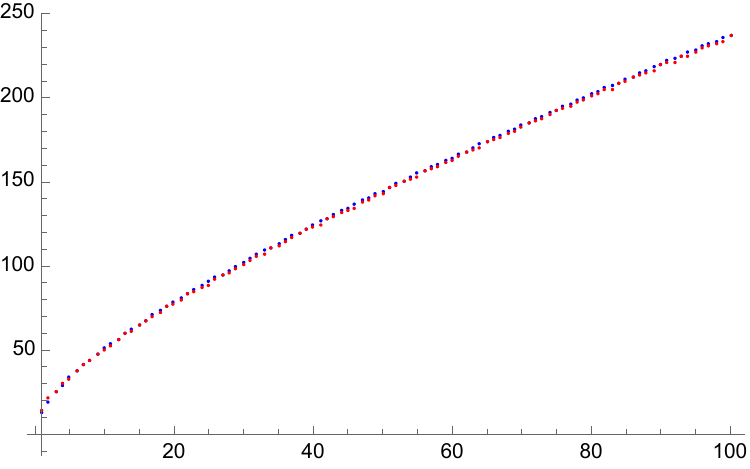}
\end{center}
\caption{The plot shows the proximity of the $n$-th zero of zeta with the $n$-th element of   Spec$D$.\label{compareplots}}
\end{figure}

The prolate spheroidal wave functions play a key role as we saw above in relation with the Riemann zeta function. In all these applications they appear as eigenfunctions of the angle operator between two orthogonal projections in the Hilbert space $L^2(\mathbb{R})^{\text {ev }}$ of even square integrable function on $\mathbb{R}$. These projections depend on a parameter $\lambda>0$, the projection $P_\lambda$ is given by the multiplication with the characteristic function of the interval $[-\lambda, \lambda] \subset \mathbb{R}$. The projection $\widehat{P_\lambda}$ is its conjugate by the Fourier transform $\mathbb{F}_{e_{\mathbb{R}}}$.
In all the above applications of prolate spheroidal wave functions the miraculous existence, discovered by the Bell Labs group  of a differential operator $PW_\lambda$ commuting with the angle operator, plays only an auxiliary role. In the present section we explain another "miracle" discovered in our joint work with H. Moscovici \cite{CM}: a careful study of the natural self-adjoint extension of $PW_\lambda$ extended  to $L^2(\mathbb{R})$ shows that it still has discrete spectrum and that its negative eigenvalues reproduce the ultraviolet behavior of the squares of zeros of the Riemann zeta function. In a similar way the positive spectrum corresponds, in the ultraviolet regime, to the trivial zeros. This coincidence holds for two values $\lambda=1$ and $\lambda=\sqrt{2}$. The conceptual reason for this coincidence is the link between the operator of \eqref{wop}, \ie
\begin{equation}
(\pw_{\lambda}\psi)(q) = \,- \partial (( \lambda^2- q^2) \partial )\,
\psi(q) + (2 \pi \lambda  q)^2 \, \psi(q)   \label{WLambdaq1}
\end{equation}
and the square of the scaling operator $\slashed S:=x \partial_x$. As we saw above the compression of the scaling $\vartheta(f)$ to Sonin's space  was shown to be the root of Weil's positivity at the archimedean place on test functions with support in the interval $\left[2^{-1 / 2}, 2^{1 / 2}\right]$, but
 	I had shown in 1998 that the prolate wave operator 
 	admits a unique selfadjoint extension commuting with the orthogonal projections $P_\lambda$ and $\widehat{P_\lambda}$. It is invariant under the Fourier transform and restricts to the Sonin space, \ie the space of square integrable even functions which vanish identically as well as their Fourier transform on the interval $[-\lambda,\lambda]$.
 	We discovered in  \cite{CM}, that this restriction of $\pw_{\lambda}$ to the Sonin space, which is a selfadjoint operator, provides a spectral realization of the ultraviolet regime of the zeros of the Riemann zeta function.	More precisely the spectrum of this operator is formed by a  discrete sequence $\nu_k$ of negative numbers and the numbers $2 \sqrt \nu_k$ for the choice $\lambda=\sqrt 2$ have the same ultraviolet behavior as the $\rho -\frac 12$ where the $\rho$'s are the zeros of zeta. \newline
 	In fact using the Darboux process we constructed a Dirac square root $D$ of $\pw_\lambda$ depending on a deformation parameter, and whose spectrum is formed of the numbers $2 \sqrt \nu_k$ for the choice $\lambda=\sqrt 2$, and has the same ultraviolet behavior as the zeros of the Riemann zeta function. Figure \ref{compareplots} shows the spectrum  of the 
 	operator $D$ and the imaginary parts of zeros of zeta. This ultraviolet spectral similarity  suggests that one has spectrally captured the contribution of the archimedean place to the mysterious zeta spectrum. The ambiguity in the choice of the Dirac square root $D$ is deeply related to the differential Galois theory of the prolate differential equation studied by J. P. Ramis and his collaborators \cite{Ramis, Ramis1}.
 	
\section{Conclusion}

The Riemann Hypothesis has catalyzed the development of vast areas of mathematics, from classical analysis to modern arithmetic geometry and mathematical physics. Each approach has enriched our understanding of the zeta function and its generalizations, even when falling short of a proof.

In this survey, we have provided a comprehensive view of this mathematical landscape. We began with a detailed examination of what is known about the Riemann zeta function and its zeros, then surveyed the remarkable variety of mathematical theories developed over 165 years in pursuit of RH—from classical analytic methods and the theory of L-functions to modern approaches via random matrix theory, operator theory, and arithmetic geometry. We explored equivalent formulations of the hypothesis, each offering its own perspective on why this problem has proven so resistant to solution. While this comprehensive survey demonstrates the richness of mathematical ideas generated by RH, it also reveals how even the most sophisticated modern approaches have yet to crack this 165-year-old puzzle.

Against this backdrop of accumulated knowledge, the second part of this paper offers something different-a return to Riemann's original viewpoint with fresh eyes. Our discovery of a large class of functions directly related to the Weil quadratic form and with zeros provably on the critical line, combined with the extraordinary numerical evidence linking truncated Euler products to the actual zeros of zeta, suggests that Riemann's original insights may contain more power than previously realized. The accuracy achieved using only primes less than 13—with errors as small as $2.6 × 10^{-55}$—cannot be dismissed as coincidence.

The geometric framework presented here, utilizing the trace formula and spectral methods, offers a potential path forward: proving that zeros of appropriately constructed approximating functions converge to zeros of zeta. This synthesis of computational discovery, classical analysis, and modern geometric techniques represents a new approach that honors Riemann's legacy while employing contemporary mathematical tools.

Whether this path leads to a proof of RH remains to be seen. However, the journey has already revealed unexpected connections between the Weil quadratic form and the world of information theory, and between computational evidence and theoretical structure. As we wrote in our letter to Riemann, sometimes the most profound truths are hidden in the simplest observations.

For comprehensive surveys of the Riemann Hypothesis and its many facets, see \cite{Patterson,BombieriRH,Edwards,IvicBook,ConreyLindelof,ConreyRH}. The present work, to be continued in collaboration with C. Consani and H. Moscovici, offers a new chapter in this ongoing story.

\begingroup
\parindent 0pt
\parskip 1ex
\theendnotes
\endgroup

\section*{Acknowledgements}
I would like to thank Jacques Dixmier and the first referee for very useful comments. 

%\newpage 

\vspace{20pt}
\scriptsize
\noindent

\textsc{Alain Connes} \\ % Small Caps for Name (optional, looks professional)
\textit{HES and College de France} \\
\textit{Paris, 75005, France} \\
\textit{E-mail address:} \texttt{alain@connes.org}

\end{document}